\documentclass[11pt]{article}
\usepackage[margin=0.86in]{geometry}
\usepackage{amsmath,amssymb,amsfonts,amsthm}
\usepackage{graphicx,float}
\usepackage[numbers]{natbib}
\usepackage[colorlinks=true,linkcolor=blue,citecolor=green,urlcolor=black,
            hypertexnames=false]{hyperref}
\usepackage{microtype}
\hypersetup{
  pdftitle={Efficient Importance Sampling of Brownian Motion Avoiding Poissonian Obstacles for Long Time},
  pdfauthor={Jose Blanchet, Paul Dupuis, Hao Liu},
  pdfkeywords={importance sampling, Wiener sausage, rare-event simulation, Poisson obstacles, asymptotic optimality}
}

\newcounter{algorithm}
\newenvironment{algorithm}[1][]{\par\addvspace{\medskipamount}\noindent\begin{minipage}{\linewidth}\refstepcounter{algorithm}}{\end{minipage}\par\addvspace{\medskipamount}}
\newcommand{\algcaption}[1]{\noindent\textbf{Algorithm \thealgorithm. #1}\par\nopagebreak[4]\smallskip\nopagebreak[4]}

\newcommand{\E}{\mathbb{E}}
\newcommand{\Pp}{\mathbb{P}}
\newcommand{\R}{\mathbb{R}}
\newcommand{\ind}{\mathbf{1}}
\newcommand{\dd}{\,\mathrm{d}}
\newcommand{\eps}{\varepsilon}
\newcommand{\BetaHat}{\widehat{\beta}}
\newcommand{\JT}{\widetilde{\mathcal J}}
\newcommand{\Jc}{\widetilde J}
\newcommand{\chic}{\widetilde\chi}
\newcommand{\vstar}{v_\star}
\newcommand{\vol}{\operatorname{Leb}_d}
\newcommand{\dist}{\operatorname{dist}}
\newcommand{\diag}{\mathrm{diag}}
\newcommand{\cem}{\mathrm{cem}}
\newcommand{\res}{\mathrm{res}}
\newcommand{\FK}{\mathrm{FK}}
\newcommand{\calN}{\mathcal N}
\newcommand{\calI}{\mathcal I}
\newcommand{\calV}{\mathcal V}
\newcommand{\calC}{\mathcal C}
\newcommand{\calE}{\mathcal E}
\newcommand{\calK}{\mathcal K}
\newcommand{\calL}{\mathcal L}
\newcommand{\esc}{\mathrm{esc}}
\newcommand{\raw}{\mathrm{raw}}
\newcommand{\sk}{\mathrm{sk}}
\newcommand{\grid}{\mathrm{grid}}
\newcommand{\up}{\mathrm{up}}
\newcommand{\cut}{\mathrm{cut}}
\theoremstyle{plain}
\newtheorem{theorem}{Theorem}[section]
\newtheorem{proposition}[theorem]{Proposition}
\newtheorem{lemma}[theorem]{Lemma}
\newtheorem{corollary}[theorem]{Corollary}
\theoremstyle{definition}
\newtheorem{definition}[theorem]{Definition}

\title{Efficient Importance Sampling of Brownian Motion Avoiding Poissonian Obstacles for Long Time}

\author{%
Jose Blanchet%
\thanks{The authors are listed in alphabetical order.}%
\textsuperscript{\phantom{*},}%
\thanks{Department of Management Science and Engineering,
Stanford University. Email: \texttt{jose.blanchet@stanford.edu}}
\and
Paul Dupuis%
\footnotemark[1]\textsuperscript{\phantom{*},}%
\thanks{Division of Applied Mathematics, Brown University.
Email: \texttt{paul\_dupuis@brown.edu}}
\and
Hao Liu%
\footnotemark[1]\textsuperscript{\phantom{*},}%
\thanks{Department of Management Science and Engineering,
Stanford University. Email: \texttt{haoliu20@stanford.edu}}
}
\date{\today}

\begin{document}
\maketitle

\begin{abstract}
We design and analyze an efficient importance sampling estimator of the
probability that, over a long time horizon, a Brownian particle avoids
obstacles of a fixed radius placed according to an independent Poisson
random field. Donsker and Varadhan identify the exponential decay rate of
this probability as the optimal balance between the volume cleared of
obstacles and the spectral cost of confinement, and the confinement
results of Sznitman and Povel show that the particle survives by
localizing its path in a ball of (roughly) deterministic radius but with
a random center. This strategy makes direct Monte Carlo exponentially
ineffective over long horizons, and makes a practical yet asymptotically
optimal sampler non-trivial to design, since the sampler must track the
random region swept out by the path. For every fixed dimension
\(d\ge2\) and every fixed obstacle radius, we construct a sequential
importance sampler that is unbiased and asymptotically optimal. The full algorithm has finite expected total work at every fixed horizon, and both its expected work and the coefficient of variation grow at most subexponentially with the time horizon.
The method avoids repeatedly solving eigenvalue
problems on the evolving explored regions and suits parallel implementation.
Numerical experiments illustrate the efficiency of the method.
\end{abstract}

\medskip
\noindent\textbf{MSC2020 subject classifications.}
Primary 65C05, 60F10; secondary 60J65, 60K37.


\tableofcontents

\section{Introduction}
\label{sec:introduction}

We study the design and analysis of an efficient importance sampling
estimator for the probability that a Brownian particle avoids Poissonian
obstacles over a long time horizon. Let $\widehat B$ be a standard Brownian
motion in a fixed dimension $d\geq2$, let $\calN$ be an independent
homogeneous Poisson point process of unit intensity on $\R^d$, and let
the obstacles be the closed Euclidean balls of a fixed radius $a > 0$ centered
at the points of $\calN$. We write $\overline B_2(x,r)$ for the closed
Euclidean ball with center $x$ and radius $r$. Conditioning on $\widehat B$,
the Poisson void formula expresses the probability of avoiding every
obstacle up to time $t$ as the exponential of minus the volume of the
region swept out by a ball of the obstacle radius travelling along the
path. That region is the Wiener sausage
\[
  W_t^r(\widehat B)
  =\bigcup_{0\leq s\leq t}\overline B_2(\widehat B_s,r).
\]

Brownian scaling converts this fixed-radius, long-horizon problem into an
equivalent shrinking-radius problem, which is the form we work with
throughout. In this rescaled formulation, we write $B$ for the standard
Brownian motion and $T$ is the time horizon after rescaling. For $a>0$, we set
\[
  \eps_T=aT^{-1/d},\qquad
  V_T=|W_T^{\eps_T}(B)|,\qquad
  \alpha_T=\mathbb E_0\big[e^{-TV_T}\big],
\]
where $\mathbb E_0$ denotes expectation for $B$.
Section~\ref{sec:setting} records the scaling identity and the
Poisson-obstacle interpretation in detail
\cite{DonskerVaradhan1975WienerSausage}.
Our goal is to design an efficient importance sampling estimator for
the quantity $\alpha_T$.

Wiener-sausage functionals and the associated geometry arise in a
wide range of applications, including, but not limited to,
\begin{itemize}
\item Molecular communication: A diffusing information molecule may be detected
by absorbing spherical receivers. When the receiver centers form a
Poisson process, the no-detection probability is a Laplace transform of
the volume swept out by the molecular path~\cite{SabuGupta2020Detection}.

\item Trapping and porous media: Static traps of different radii,
intermittent traps arising in models of T-cell selection, and finite-sized Brownian particles in
porous media lead to related survival and encounter models
\cite{BerezhkovskiiMakhnovskiiWeiss2001Radii,Kosmrlj2011Thymic,
Torquato1991FiniteParticles}.

\item Chemical, photophysical, and biological reactions:
Wiener-sausage constructions are used to estimate particle-collision
probabilities, quantify the volume explored by excitons in polymer-dot
photocatalysts, and simulate encounters of diffusing molecules with immobile
binding sites
\cite{PerfileevLyakishev2025Collision,BrnovicEtAl2025Exciton,
SigautEtAl2014Messages}.

\item Polymers and random media: Range-penalized random walks are discrete
analogues of the same volume penalty and admit both poor-solvent polymer and
annealed obstacle-survival interpretations, with the penalty favoring compact
paths \cite{Bouchot2024Range,BerestyckiCerf2021Range}.
\end{itemize}

The Wiener sausage has been studied extensively in probability theory.
Foundational work connected its volume with electrostatic capacity and
developed the theory of its fluctuations and self-intersections
\cite{Spitzer1964Capacity,LeGall1988Fluctuation,
LeGall1990SelfIntersection}.
Donsker and Varadhan established the leading Laplace asymptotics through a
variational balance between the volume penalty and the spectral cost of
Brownian confinement
\cite{DonskerVaradhan1975WienerSausage}.
Small- and shrinking-radius asymptotics place our shrinking-radius setting in
a broader classical context
\cite{Sznitman1987SmallRadius,Sznitman1990Shrinking}.
In the equivalent Poisson-obstacle formulation, Lifshitz-tail and spectral
methods characterize the survival scale, while confinement results describe
the obstacle-free region predicted by this balance and the compact geometry
of paths conditioned on long survival
\cite{Sznitman1990LifschitzI,Sznitman1990LifschitzII,
Sznitman1991Confinement,Povel1999Confinement,
Sznitman1998Obstacles,Fukushima2008Sausage}.
Lower-volume deviations and intersections of independent sausages further
describe the rare geometry relevant to one- and two-path estimates
\cite{VanDenBergBolthausenDenHollander2001Moderate,
VanDenBergBolthausenDenHollander2004Intersection}.
More recent work has investigated torsional rigidity of complements generated
by Brownian paths and shrinking sausages
\cite{VanDenBergBolthausenDenHollander2018Torsional},
the capacity of the four-dimensional Wiener sausage
\cite{AsselahSchapiraSousi2019Capacity},
long-time sausage-volume behavior for diffusions on Dirichlet spaces
\cite{Okamura2020Dirichlet},
spectral gaps in Poissonian obstacle media
\cite{Sznitman2023SpectralGap},
and connectivity properties of vacant sets generated by Wiener sausages and
Brownian interlacements
\cite{MuSapozhnikov2024Vacant}.

The existing literature is primarily concerned with asymptotic and structural
properties of the Wiener sausage and the associated survival problem. 
In
particular, it identifies the exponential scale of the survival probability
and the compact geometry of the paths that dominate it. More precisely, in the scaled variables introduced above,
\(\alpha_T=\exp\{-JT+o(T)\}\) for a constant \(J>0\), and conditional
on survival the path is typically confined to a ball with an
asymptotically deterministic radius and a random center. It does not, however,
provide a finite-horizon construction for sampling such paths together with an
efficient Monte Carlo estimator of that probability. This leads naturally
to the computational problem studied here. Direct Monte Carlo simulation is
ineffective in the rare regime: typical paths explore too much volume and
therefore contribute values of $e^{-T V_T}$ that are negligible relative to
$\alpha_T$, while the expectation is governed by atypically compact paths. 

Our goal is therefore to construct an unbiased importance-sampling estimator
of \(\alpha_T\) that is asymptotically optimal (also called logarithmically
efficient in the rare-event simulation literature)
\cite{Siegmund1976SequentialTests,GlynnWhitt1992Efficiency,
BlanchetLam2012StateDependent}. Roughly speaking, asymptotic optimality
requires the estimator's second moment to decay at the same exponential rate
as \(\alpha_T^2\) on the scale \(T\).
Definition~\ref{def:log-eff} gives the precise criterion used here.

Importance sampling for rare-event simulation has been developed through
several complementary principles. Foundational work on stochastic simulation
and asymptotic efficiency established the likelihood-ratio framework and
emphasized that performance must ultimately be assessed through variance or
second-moment control
\cite{GlynnIglehart1989StochasticSimulations,
GlynnWhitt1992Efficiency}.
Large-deviation analysis provides a classical route to constructing
asymptotically efficient changes of measure
\cite{SadowskyBucklew1990LargeDeviations}.
Adaptive, sequential, and state-dependent schemes instead allow the change of
measure to evolve with the simulated state
\cite{BaggerlyEtAl1999Adaptive,BlanchetLam2012StateDependent}.
Variational representations for exponential Brownian functionals, together
with stochastic-control, differential-game, and subsolution formulations,
provide related dynamic foundations and constructive methods
\cite{BoueDupuis1998Variational,DupuisWang2004DifferentialGames,
DupuisWang2005Dynamic,DupuisWang2007Subsolutions,
VandenEijndenWeare2012SmallNoise}.
These positive results are accompanied by an important caution: matching only
the leading static rare-event mechanism need not control the
likelihood-weighted second moment and can lead to poor or even infinite
variance
\cite{GlassermanWang1997Counterexamples}.

For this functional, the natural sequential benchmark is
the exact zero-variance law. Given the path history at an intermediate stage,
let $H_T(S_k)$ be the conditional expectation of the remaining exponential functional value. If $F_k$ is the exact incremental exponential functional value over the next stage and
$S_{k+1}$ is the resulting state, then the strong Markov property gives the
Bellman identity
\[
  H_T(S_k)
  =
  \mathbb{E}\!\left[
    F_k H_T(S_{k+1})
    \,\middle|\,
    S_k
  \right].
\]
The ideal transition biases each possible continuation proportionally to
$F_k H_T(S_{k+1})$. This is the path-dependent counterpart of the
zero-variance principle for Markov-process expectations
\cite{AwadGlynnRubinstein2013ZeroVariance}
and is closely related to conditioned-process and generalized Doob-transform
constructions
\cite{ChetriteTouchette2015Conditioned}.

The large-deviation rate function \cite{DonskerVaradhan1975WienerSausage} identifies a tradeoff between the volume explored by the path and the spectral cost of confining Brownian motion in the explored region.
The volume--confinement picture therefore suggests approximating the
continuation value through principal spectral objects: on a fixed domain, a
principal Dirichlet eigenfunction and its Doob transform bias Brownian motion
toward long survival. The enlarged explored region is also recorded sequentially. For the Wiener sausage, however, the relevant state
includes the irregular region explored so far. After each enlargement, an
adaptive spectral implementation would have to update a path-dependent
boundary continuation value, principal eigenvalue and eigenfunction, and a
local Doob transform on the new geometry. Repeatedly solving these spectral problems on the evolving explored
region is computationally demanding.
Moreover, establishing asymptotic optimality for a numerical plug-in law would
require controlling the accumulated approximation error in the
estimator's second moment, which is difficult
\cite{LEcuyerTuffin2011ZeroVariance}.

Our main contribution is to construct an unbiased and asymptotically optimal importance sampling method that does not rely on the numerical approximation of the exact
continuation value $H_T$. Instead, we divide the time horizon into fine time steps and group them into coarser blocks. At each block, we augment the transition
with a finite cloud of conditionally independent candidate Brownian
path segments sampled on the fine time grid. Each candidate is scored by its immediate exponential new-volume
weight multiplied by a positive continuation guide defined
in~\eqref{eq:guide}, and one candidate is selected with probability
proportional to this score. Crucially,
the empirical mean of the candidate scores is retained in the estimator
weight. In the terminology of sequential Monte Carlo, the resulting weighted
candidate is properly weighted \cite{LiuChen1998SequentialMonteCarlo}. The
random empirical normalizer thus yields an exact one-step proper-weight
identity. The finite-cloud step is related in spirit to sequential particle
and twisted-particle constructions
\cite{DelMoralGarnier2005Genealogical,WhiteleyLee2014Twisted}, but here a
single selected lineage is followed and the empirical normalizer is retained
to obtain an exact finite-horizon identity. The guide is not used as an uncontrolled
plug-in approximation to the ideal population normalizer. It induces an implicit change of measure through weighted candidate selection and is designed to control the estimator’s second moment. We also show that
the full algorithm has finite expected total work at every fixed
horizon, and this expectation grows at most subexponentially
in \(T\) for fixed \(d\) and \(a\).

A close problem-specific comparison is the lattice static-trap study, which uses a compactness-biased
whole-path polymer Monte Carlo method followed by reweighting
\cite{BarkemaBiswasVanBeijeren2001StaticTraps}.
A separate direct numerical study uses Brownian-dynamics simulation to investigate the statistical
distribution of Wiener-sausage volume
\cite{YangMakhnovskiiSheuLin2000Simulation}. Such a study describes typical
volume fluctuations, whereas the Laplace transform considered here is governed
by the far lower tail; accurate simulation of the former therefore does not by
itself provide an efficient estimator of the latter.
These works do not provide the
unbiasedness and asymptotic optimality developed here.

Our next contribution is the numerical validation of the proposed construction. We implement a finite-horizon lattice surrogate in which the continuous Brownian path is replaced by a nearest-neighbor random walk on $\mathbb{Z}^d$, and the Wiener sausage is replaced by the exact union of discrete lattice balls centered at the sites visited by the walk.
 The experiments report the estimation of very small surrogate expectations and illustrate the practical parallelism enabled by batched geometry. The results show substantial variance reduction and computational gains relative to direct Monte Carlo. 

The rest of the paper is organized as follows. Section~\ref{sec:model} defines the setting, derives the zero-variance Bellman law, and describes the natural
spectral/Doob-transform approach. Section~\ref{sec:algorithm} specifies our
method. Sections~\ref{sec:exactness} --
\ref{sec:efficiency} establish exact finite-horizon unbiasedness, derive an
exact finite-$T$ second-moment representation, and prove the asymptotic
optimality of our method. Section~\ref{sec:numerics} reports the
lattice-surrogate experiments, and
Section~\ref{sec:discussion} concludes. The appendices
contain the technical details of the construction and the asymptotic analysis and additional details about the experiments.

\section{The target functional and the ideal change of measure}
\label{sec:model}

\subsection{Brownian model and natural-radius target}\label{sec:setting}

Fix \(d\ge2\). On
\(\Omega_T=C_0([0,T],\R^d)\) equipped with the uniform topology and Borel sigma algebra, let \(B_t(\omega)=\omega(t)\), let
\(\mathcal F_t=\sigma(B_s:0\le s\le t)\), and let \(\Pp_T\) be Wiener
measure. Thus \(B_0=0\), \(\operatorname{Cov}_{\Pp_T}(B_t)=tI_d\), and
the generator is \(\Delta/2\). For \(r>0\) and \(0\le s\le t\le T\),
define the closed-ball Wiener sausage
\begin{equation*}
  W_{[s,t]}^r(B)
  =\bigcup_{s\le u\le t}\overline B_2(B_u,r),
  \qquad W_t^r=W_{[0,t]}^r.
\end{equation*}
The range of a continuous path on a compact interval is compact, so the
sausage is compact and its volume is a measurable path functional.

Fix also \(a>0\), and for every \(T>0\) put
\begin{equation}
  \eps_T=aT^{-1/d},\qquad
  V_T=\vol(W_T^{\eps_T}),\qquad
  F_T=e^{-TV_T},\qquad
  \alpha_T=\E_{\Pp_T}F_T,
  \label{eq:target}
\end{equation}
where for any Lebesgue measurable set $A \subset \mathbb{R}^d$, $\vol(A)$ is its Lebesgue measure in $\mathbb{R}^d$. 

Brownian scaling converts \eqref{eq:target} into a fixed-radius, long-time
problem. Indeed, if \(\widehat T=T^{1+2/d}\), then
\(\widehat B_s=T^{1/d}B_{T^{-2/d}s}\), \(0\le s\le\widehat T\), is a
standard Brownian motion and
\begin{equation}
  T\vol(W_T^{\eps_T}(B))
  =\vol(W_{\widehat T}^{a}(\widehat B)),
  \qquad
    \alpha_T=\E\exp\{-\vol(W_{\widehat T}^{a}(\widehat B))\}.
  \label{eq:scaling}
\end{equation}
The volume identity in~\eqref{eq:scaling} holds pathwise by the
definition of \(\widehat B\). This fixed-radius functional has the standard Poisson-obstacle
interpretation. Recall that \(\mathcal N\) is a unit-intensity homogeneous Poisson
point process on \(\R^d\), independent of \(\widehat B\), and let
\[
  \mathcal O^a
  =\bigcup_{z\in\mathcal N}\overline B_2(z,a)
\]
be the corresponding obstacle set. Under the joint law
\(\mathbb P^{\mathrm{ann}}\), survival through time
\(\widehat T\) is the event that \(\mathcal N\) has no point in
\(W_{\widehat T}^a(\widehat B)\). Hence the Poisson void formula gives
\[
  \mathbb P^{\mathrm{ann}}\!\left(
    \widehat B_s\notin\mathcal O^a
    \text{ for all }0\le s\le\widehat T
    \,\middle|\,\widehat B
  \right)
  =e^{-\vol(W_{\widehat T}^{a}(\widehat B))}.
\]
Thus \(\alpha_T\) is the annealed survival probability of a Brownian particle
among fixed-radius Poissonian obstacles, written in natural-radius variables
\cite{DonskerVaradhan1975WienerSausage,Sznitman1998Obstacles}.

Let \(\lambda_d\) be the principal Dirichlet eigenvalue of \(-\Delta\)
on a Euclidean ball of volume one. Since the Brownian generator is
\(\Delta/2\), define
\begin{equation*}
  \chi_d=\frac{\lambda_d}{2},
  \qquad
  \Lambda_{\FK}(v)=\chi_dv^{-2/d}.
\end{equation*}
This factor of \(1/2\) is used throughout. We also define the constants
\begin{equation}
  \vstar
  =
  \left(\frac{\lambda_d}{d}\right)^{d/(d+2)},
  \qquad
  J
  =
  \inf_{v>0}\bigl\{v+\Lambda_{\FK}(v)\bigr\}
  =
  \frac{d+2}{2}\vstar,
  \qquad
  p=\frac{d}{d+2}.
  \label{eq:rate-constants}
\end{equation}

The classical fixed-radius Donsker--Varadhan Laplace asymptotic,
combined with the Brownian scaling in \eqref{eq:scaling}, gives the
following natural-radius reformulation
\cite{DonskerVaradhan1975WienerSausage}.

\begin{proposition}[Natural-radius target rate]
\label{prop:target-rate}
For every fixed \(d\ge2\) and \(a>0\),
\begin{equation}
  \log\alpha_T=-JT+o(T),
  \qquad T\to\infty.
  \label{eq:target-rate}
\end{equation}
\end{proposition}

The variational expression \(J\) has the classical
volume--confinement interpretation. In the natural-radius formulation
(more generally, the Poisson-obstacle interpretation holds for any deterministic
\(\eps_T>0\), while \(\eps_T=aT^{-1/d}\) is the choice that
rescales to the classical fixed-radius setting),
the factor $e^{-T|W_T^{\eps_T}|}$
can be interpreted as a Poisson void probability with intensity \(T\):
conditional on the Brownian path, it is the probability that no Poisson
obstacle center lies in the Wiener sausage. Thus, if the favorable
obstacle-free region has volume approximately \(v\), its leading
Poisson clearing cost is $ \exp\{-Tv+o(T)\}.$
Hence the first term \(v\) in $ v+\Lambda_{\FK}(v)$
represents the obstacle-free volume cost on the exponential scale \(T\).

The second term represents the cost of confining the Brownian path
inside such a region over the time horizon \(T\). Among domains of
volume \(v\), the Faber--Krahn inequality shows that a Euclidean ball
minimizes the principal Dirichlet eigenvalue. Since the Brownian
generator is \(\Delta/2\), the least confinement rate is $\chi_d v^{-2/d}
  =
  \Lambda_{\FK}(v).$
The corresponding leading confinement cost is therefore $\exp\!\left\{
    -T\Lambda_{\FK}(v)
  \right\}.$

Thus the two terms in the variational expression have the
interpretation
\[
  \underbrace{v}_{\text{obstacle-free volume cost}}
  +
  \underbrace{\Lambda_{\FK}(v)}_{\text{spectral confinement cost}},
\]
and together they give the leading exponential cost.
Optimizing this balance over \(v>0\) gives the unique minimizer \(\vstar\)
and the value \(J\) in \eqref{eq:rate-constants}. In this way,
Proposition~\ref{prop:target-rate} states that the leading decay of the target is governed by the optimal
balance between obstacle-free volume and Brownian confinement.
This interpretation also motivates the structure of the importance-sampling
construction below: at an intermediate time, one must account for both the
volume already explored and the spectral cost of confining the path over the
remaining horizon. 

\subsection{Continuation value and zero-variance benchmark}
\label{sec:ideal-law}

We next describe the exact sequential law that would solve the simulation
problem if its continuation values were available. At time \(t\), retain the
complete observed prefix \(\omega_{[0,t]}\), its endpoint \(x=\omega(t)\), and
a certified Borel set \(D\subset W_t^{\eps_T}(\omega)\) whose volume has
already been included in the exponential weight (paid inner region). Write $S=(t,\omega_{[0,t]},D).$
The complete prefix is needed because, when \(D\) is only an inward paid
subset, the reduced information \((t,x,D)\) does not determine the uncharged
part of the past continuous sausage. This can also be viewed as the use of explored regions as an inner approximation of the Wiener sausage.

Starting from state \(S\), the remaining exponential functional is
\[
  \exp\{-T(\vol(W_T^{\eps_T})-\vol(D))\}.
\]
Its conditional expectation defines the exact continuation value
\begin{equation*}
  H_T^\star(S)
  =\E_{\Pp_T}\!\left[
      \exp\{-T(\vol(W_T^{\eps_T})-\vol(D))\}
      \,\middle|\,\mathcal F_t
    \right](\omega).
\end{equation*}
At the initial state \(S_0=(0,(0),\varnothing)\), the continuation value
reduces to the target: \(H_T^\star(S_0)=\alpha_T\).

Let \(\Gamma\) denote the next Brownian segment, stopped no later than the remaining horizon, and let $S^+=\Phi(S,\Gamma)$
be the successor state obtained by appending \(\Gamma\) to the observed
prefix and applying a measurable, adapted update of the paid set.
Write \(D^+(\Gamma)\supset D\) for the updated certified subset of the
successor sausage. The incremental exponential weight is
\[
  F_{\rm inc}(D,\Gamma)
  =\exp\{-T(\vol(D^+(\Gamma))-\vol(D))\}.
\]
The remaining functional factors into the current incremental weight and the
functional remaining from \(S^+\). Therefore, the tower property and the
independence of Brownian increments give the Bellman identity
\begin{equation}
  H_T^\star(S)
  =\E_{\Pp_T}\!\left[
       F_{\rm inc}(D,\Gamma)H_T^\star(S^+)
       \,\middle|\,S
    \right].
  \label{eq:bellman}
\end{equation}
In particular, the right-hand side is the normalizing constant of the ideal
one-stage transition. If \(\Pp_T(\dd\Gamma\mid S)\) denotes the ordinary
Brownian continuation kernel, define
\begin{equation}
\begin{aligned}
  \frac{\dd\mathbb Q_T^\star(\cdot\mid S)}
       {\dd\Pp_T(\cdot\mid S)}(\Gamma)
  &=\frac{F_{\rm inc}(D,\Gamma)H_T^\star(\Phi(S,\Gamma))}
    {\displaystyle\int F_{\rm inc}(D,\gamma)
      H_T^\star(\Phi(S,\gamma))\,\Pp_T(\dd\gamma\mid S)}=\frac{F_{\rm inc}(D,\Gamma)H_T^\star(S^+)}{H_T^\star(S)}.
\end{aligned}
  \label{eq:zero-variance}
\end{equation}
The denominator is the ordinary Brownian expectation of the complete
one-stage score; \eqref{eq:bellman} identifies it with
\(H_T^\star(S)>0\). Thus \eqref{eq:zero-variance} integrates to one. It
biases each possible continuation in proportion to its current exponential
weight times its exact future value, which is the standard zero-variance
principle for Markov-process expectations
\cite{AwadGlynnRubinstein2013ZeroVariance}.

To see the zero-variance property, iterate \eqref{eq:zero-variance} along
states \(S_0,S_1,\ldots,S_K\), and write
\(F_{{\rm inc},k}=F_{\rm inc}(D_k,\Gamma_k)\). The
reference-over-proposal likelihood ratio satisfies
\[
  \prod_{k=0}^{K-1}
  \frac{\dd\Pp_T(\cdot\mid S_k)}
       {\dd\mathbb Q_T^\star(\cdot\mid S_k)}(\Gamma_k)
  =
  \frac{H_T^\star(S_0)}
       {\left(\prod_{k=0}^{K-1}F_{{\rm inc},k}\right)
        H_T^\star(S_K)}.
\]
At terminal time,
\[
  \left(\prod_{k=0}^{K-1}F_{{\rm inc},k}\right)
  H_T^\star(S_K)
  =F_T,
\]
because the product charges \(\vol(D_K)\) and the terminal continuation value
contains exactly the remaining volume
\(\vol(W_T^{\eps_T})-\vol(D_K)\). Hence the likelihood-corrected terminal
functional is the constant
\[
  F_T
  \prod_{k=0}^{K-1}
  \frac{\dd\Pp_T(\cdot\mid S_k)}
       {\dd\mathbb Q_T^\star(\cdot\mid S_k)}(\Gamma_k)
  =H_T^\star(S_0)=\alpha_T.
\]
This law is therefore a benchmark rather than an implementation, since evaluating
\(H_T^\star\) already includes the original problem at \(S_0\).

\subsection{The spectral route and its adaptive obstruction}
\label{sec:spectral}

\subsubsection{A principal-eigenfunction benchmark}
The volume--confinement balance derived above suggests a natural idea of using
spectral information to represent or approximate future value. To separate
the two ingredients of such a construction, first consider a bounded regular
domain \(U\subset\R^d\). Let \(\lambda_{1/2}(U)\) be the principal Dirichlet
eigenvalue of \(-\Delta/2\), and let \(\varphi_U\) be a corresponding positive
principal eigenfunction:
\[
  -\frac12\Delta\varphi_U=\lambda_{1/2}(U)\varphi_U
  \quad\text{in }U,
  \qquad
  \varphi_U=0\quad\text{on }\partial U.
\]
Thus \(\lambda_{1/2}(U)\) is one half of the principal eigenvalue of
\(-\Delta\) on \(U\). It determines the exponential cost of keeping
Brownian motion inside \(U\), while \(\varphi_U\) gives the spatial bias of the
corresponding long-survival Doob transform. Thus the eigenpair describes
confinement in a fixed region.

To account for exit and subsequent evolution, let \(\mathbb P_x\) denote Wiener
measure on \(C([0,\infty),\mathbb R^d)\) for Brownian motion $B$, and let
\(\mathbb E_x\) denote the corresponding expectation. Define the stopping time $\tau_U=\inf\{r\ge0:B_r\notin U\}.$ For \(\theta<\lambda_{1/2}(U)\) and a bounded positive continuous
boundary payoff \(g\) on \(\partial U\), define
\[
  h_U^g(x;\theta)
  =\E_x\!\left[e^{\theta\tau_U}g(B_{\tau_U})\right].
\]Then under standard boundary regularity,
\[
  \left(\frac12\Delta+\theta\right)h_U^g=0
  \quad\text{in }U,
  \qquad h_U^g=g\quad\text{on }\partial U,
\]
and the corresponding Doob transform has drift
\(\nabla\log h_U^g\). The factor \(e^{\theta\tau_U}\) compensates for a
candidate exponential confinement rate, whereas \(g\) determines how the law
values different exit locations. In a finite-horizon problem the relevant boundary datum generally
also depends on the exit time, so one needs a family \(g(r,y)\), not only a
single function of \(y\).

For the Wiener-sausage problem, if a local segment \(\Gamma\) exits \(U\)
after time \(r\) at location \(y\), then its exact one-segment score is
\begin{equation*}
  \mathfrak g_T^\star(S;\Gamma)
  =F_{\rm inc}(D,\Gamma)H_T^\star(S^+).
\end{equation*}
Accordingly, the oracle boundary datum would have to be a regular conditional
version of
\[
  g_{T,S}^\star(r,y)
  =e^{-\theta r}
   \E_{\Pp_T}\!\left[
     \mathfrak g_T^\star(S;\Gamma)
     \,\middle|\,
     \tau_U=r,\ B_{\tau_U}=y,\ S
   \right].
\]
The factor \(e^{-\theta r}\) removes the factor already included in
\(e^{\theta\tau_U}\). This formula makes the division of labor explicit:
the principal spectral data quantify the cost and spatial profile of
confinement in the current region, while the boundary datum determines when
and where an exit is valuable by incorporating the new volume created by the
segment and the continuation value of the future enlarged geometry.
This is the local spectral/Doob-transform counterpart of the exact sequential
law in Section~\ref{sec:ideal-law}.

\subsubsection{Why a plug-in spectral law is difficult here}
The goal is to design an importance sampling estimator with controlled normalized second moment. For any unbiased importance-sampling estimator \(Z_T\) under a law
\(Q_T\), the normalized second moment satisfies
\[
  \frac{\E_{Q_T}Z_T^2}{\alpha_T^2}
  =1+\frac{\operatorname{Var}_{Q_T}(Z_T)}{\alpha_T^2}\ge1.
\]
To see its connection with
sampling effort, let \(\overline Z_{T,n}\) be the average of \(n\)
i.i.d. copies of \(Z_T\). Then its relative mean squared error is
\[
  \frac{\E[(\overline Z_{T,n}-\alpha_T)^2]}{\alpha_T^2}
  =
  \frac1n\left(
    \frac{\E_{Q_T}Z_T^2}{\alpha_T^2}-1
  \right).
\]
The relative variance $\frac{\operatorname{Var}_{Q_T}(Z_T)}{\alpha_T^2}$ therefore determines the number of independent
samples needed to attain a prescribed relative mean squared error.
We seek a sampling law for which this sample requirement grows at most
subexponentially in \(T\). The asymptotic-optimality goal, formalized in
Definition~\ref{def:log-eff}, is
\[
  \frac1T\log\frac{\E_{Q_T}Z_T^2}{\alpha_T^2}\longrightarrow0,
  \qquad T\to\infty.
\]
Since \(\alpha_T^2=e^{-2JT+o(T)}\), this condition is equivalent to
\(\E_{Q_T}Z_T^2=e^{-2JT+o(T)}\): the second moment must decay at the
same exponential rate as the squared target.
For a plug-in spectral construction, we
therefore need to understand how approximating the continuation value
affects these weights and their second moment.

A direct adaptive implementation would replace \(H_T^\star\) by a positive
finite measurable function \(\widehat H_T\). For a finite sequence of stages ending at
\(T\), write \(F_k=F_{\rm inc}(D_k,\Gamma_k)\),
\(S_{k+1}=\Phi(S_k,\Gamma_k)\), and
\[
  (\mathcal B_{T,k}f)(S_k)
  =\E_{\Pp_{T,k}}\!\left[
     F_kf(S_{k+1})\,\middle|\,S_k
   \right],
\]
where \(\Pp_{T,k}(\cdot\mid S_k)\) denotes the conditional law of
the next Brownian path segment given \(S_k\). Assume the normalizers
\((\mathcal B_{T,k}\widehat H_T)(S_k)\) are finite, then the normalized
plug-in transition is
\[
  \frac{\dd\widehat{\mathbb Q}_{T,k}(\cdot\mid S_k)}
       {\dd\Pp_{T,k}(\cdot\mid S_k)}
  =\frac{F_k\widehat H_T(S_{k+1})}
         {(\mathcal B_{T,k}\widehat H_T)(S_k)}.
\]
Let \(\widehat{\mathbb Q}_T\) be the resulting path law and define its
logarithmic Bellman residual by
\[
  R_{T,k}(S_k)
  =\log\frac{(\mathcal B_{T,k}\widehat H_T)(S_k)}
                  {\widehat H_T(S_k)}.
\]
The likelihood-corrected estimator satisfies the exact telescoping identity
\[
\begin{aligned}
  \widehat Z_T^{\rm plug}
  &:=F_T\prod_{k=0}^{K-1}
      \frac{\dd\Pp_{T,k}(\cdot\mid S_k)}
           {\dd\widehat{\mathbb Q}_{T,k}(\cdot\mid S_k)}(\Gamma_k)=\widehat H_T(S_0)
    \frac{H_T^\star(S_K)}{\widehat H_T(S_K)}
    \exp\!\left\{\sum_{k=0}^{K-1}R_{T,k}(S_k)\right\}.
\end{aligned}
\]
Here the terminal ratio remains because
\(F_T=(\prod_{k=0}^{K-1}F_k)H_T^\star(S_K)\).
Consequently,
\[
\begin{aligned}
  \frac1T\log\frac{\E_{\widehat{\mathbb Q}_T}
       [(\widehat Z_T^{\rm plug})^2]}{\alpha_T^2}
  &=\frac2T\log\frac{\widehat H_T(S_0)}{\alpha_T}+\frac1T\log\E_{\widehat{\mathbb Q}_T}\!\left[
     \left(\frac{H_T^\star(S_K)}{\widehat H_T(S_K)}\right)^2
     \exp\!\left\{2\sum_{k=0}^{K-1}R_{T,k}(S_k)\right\}
   \right].
\end{aligned}
\]
Thus matching the initial exponential scale of \(\widehat H_T\) does not
suffice: one must also control the joint contribution of the accumulated
Bellman residual and the terminal ratio under the adaptive proposal.

This requirement is especially difficult for the Wiener sausage. The state
contains the current position, the remaining time, the full relevant path
history, and an evolving irregular explored region. After each selected
segment, a direct spectral implementation based on the updated geometry would
have to recompute a principal spectral object, a time- and location-dependent
boundary continuation datum, and the resulting local Doob transition. These
solves are sequentially coupled because the next geometry depends on the
previously selected path. They are consequently computationally burdensome,
poorly aligned with parallel candidate generation, and difficult to analyze
at the level of the accumulated second moment. Even if spectral information were used only as a guide within the candidate-augmented scheme as we discussed later, each selected segment would still require updating spectral and boundary-continuation objects on the evolving irregular region. Thus augmentation could avoid plug-in Bellman-residual control, but not the sequential and computationally burdensome spectral solves.

This is a setting-specific obstruction, not a claim that spectral importance
sampling in general lacks asymptotic-optimality guarantees. Our augmentation method below
does not attempt to compute \(H_T^\star\), a plug-in approximation
\(\widehat H_T\), or its population normalizer. Instead, it draws a finite
cloud of ordinary Brownian continuations, scores them by
\[
  \text{incremental exponential weight}
  \times
  \text{positive continuation guide},
\]
selects one candidate proportionally to its score, and retains the random
empirical normalizer in the exact proper weight. The guide steers candidate selection, while the augmented finite-cloud
normalization yields an exact one-step proper-weight identity.
The guide, grid parameters, and candidate counts are chosen to
control the second moment. Our method avoids the need to control
the accumulated residual arising from a plug-in approximation
of the exact continuation value.

\section{Sequential candidate-augmented importance sampling}
\label{sec:algorithm}

\subsection{Deterministic parameters and inward grid}

The finite law is defined for every \(T>0\).  Set
\(T_{\rm reg}=1\vee T\).  We choose approximations of the scalar ball-based variational/spectral constants via positive constants
\(\widetilde v_{\star,T},\chic_{d,T},\Jc_T>0\) which satisfy
\begin{equation}
  \Jc_T=\frac{d+2}{2}\widetilde v_{\star,T},
  \qquad
  \chic_{d,T}=\frac d2\widetilde v_{\star,T}^{\,1+2/d},
  \qquad
  |\widetilde v_{\star,T}-\vstar|
  +|\chic_{d,T}-\chi_d|+|\Jc_T-J|
  \le T_{\rm reg}^{-4}.
  \label{eq:certified-constants}
\end{equation}
Appendix~\ref{app:exact-real} gives a deterministic construction of these
constants. Only spectral data for a Euclidean ball are required and can be
certified once. By contrast, the hypothetical adaptive benchmark in
Section~\ref{sec:spectral} would require spectral calculations on evolving
irregular geometries, not only for the known Euclidean balls.

We begin by defining the fine time mesh, the coarser sequential block time
schedule, and the finite spatial grid used by the method. Put
\begin{equation}
  \ell_T=1\vee\log T,\qquad
  \delta_T=\frac{\eps_T}{8d},\qquad
  N_T=\left\lceil T^2\ell_T^3\delta_T^{-2}\right\rceil,\qquad
  h_T=\frac{T}{N_T}.
  \label{eq:raster-schedule}
\end{equation}
Here \(\ell_T\) provides logarithmic slack in the discretization,
while \(\delta_T\) is the side length of each spatial grid cell.
The interval \([0,T]\) is divided into \(N_T\) fine time gaps, each of
duration \(h_T\). Thus \(N_T\) is the number of Brownian increments and \(h_T\) is the time resolution of the selected
skeleton. The choice of \(N_T\) makes the time mesh much finer than the
Brownian time scale associated with one spatial cell; indeed, $h_T
  \le
  \frac{\delta_T^2}{T\ell_T^3}.$
This fine resolution is used later to control the continuous Brownian
motion between successive times.

To group these fine increments into sequential importance-sampling
stages, define the dimension-dependent analytic scale
\begin{equation}
\begin{aligned}
  \mathfrak H_T&=
  \begin{cases}
    \max\{\delta_T^{-(d-2)},\log N_T,1\},&d\ge3,\\
    \max\{\log(T\delta_T^{-2}),\log N_T,1\},&d=2,
  \end{cases}
  &\qquad
  b_T&=\sqrt{T\mathfrak H_T}.
\end{aligned}
\label{eq:block-rate}
\end{equation}
The quantity \(b_T\) is the nominal block duration. It is chosen
as an intermediate scale: asymptotically, each block is long enough
for the required lower-volume estimates, while \(b_T=o(T)\), so the
path is still divided into an increasing number of sequential stages.

Balance the \(N_T\) gaps among these stages by
\begin{equation*}
\begin{aligned}
  K_T&=\min\{N_T,\lceil T/b_T\rceil\},&
  j_k&=\left\lfloor\frac{kN_T}{K_T}\right\rfloor,&
  t_k&=\frac{Tj_k}{N_T},\\
  n_{T,k}&=j_{k+1}-j_k,&
  b_{T,k}&=t_{k+1}-t_k,&
  q_k&=1-\frac{t_k}{T}.
\end{aligned}
\end{equation*}
Here \(K_T\) is the number of sequential blocks, \(j_k\) is the fine
grid index at which block \(k\) begins, and \(t_k\) is its physical
starting time. The block contains \(n_{T,k}\) fine Brownian increments
and has actual duration \(b_{T,k}\). Finally, \(q_k\) is the fraction
of the total time horizon remaining at the beginning of block \(k\) which will enter the continuation guide. The floor convention distributes
the \(N_T\) grid gaps as evenly as possible across the \(K_T\)
blocks. In particular, \(j_0=0\), \(j_{K_T}=N_T\), all block durations
are positive, and every grid gap belongs to exactly one block. 

We next define the finite spatial window and the tokenization of grid
points. For
\[
  L_T^{\rm win}
  =
  1\vee\left\lfloor T^2/\delta_T\right\rfloor,
  \qquad
  \calC_T
  =
  \{-L_T^{\rm win},\ldots,L_T^{\rm win}-1\}^d,
\]
let
\[
  Q_{T,m}
  =
  \delta_T(m+[0,1)^d),
  \qquad
  \calV_T
  =
  \bigcup_{m\in\calC_T}Q_{T,m}.
\]
Thus \(\calC_T\) is the finite set of valid cell labels,
\(Q_{T,m}\) is the half-open grid cell carrying label \(m\), and
\(\calV_T\) is the corresponding spatial window. Its coordinate
half-width is approximately \(T^2\). The Borel token map is
\begin{equation*}
  \tau_T(x)=
  \begin{cases}
    \lfloor x/\delta_T\rfloor,&
      \lfloor x/\delta_T\rfloor\in\calC_T,\\
    \dagger,&\text{otherwise},
  \end{cases}
\end{equation*}
This map sends a continuous state \(x\) to a token (a discrete label identifying a grid cell).
A point inside \(\calV_T\) receives the label of its grid cell,
whereas \(\dagger\) records that the point lies outside the finite
window. The symbol \(\dagger\) is only an outside-window flag (referred to below as the dagger token).

Define $d_T^{\cem}(x)
  =
  \dist_\infty(x,\calV_T^c)$
for \(x\in\calV_T\), and set it equal to zero otherwise. This quantity
measures the distance from \(x\) to the boundary of the token window
and will be used in the positive escape component of the continuation
guide.

The spatial mesh is chosen sufficiently fine to provide a pathwise
inner certificate. Indeed, because
\(\delta_T\sqrt d<\eps_T\), if \(x\in Q_{T,m}\), then every point of
\(\overline Q_{T,m}\) is within Euclidean distance
\(\delta_T\sqrt d\) of \(x\), and hence
\begin{equation}
  \overline Q_{T,m}
  \subset
  \overline B_2(x,\eps_T).
  \label{eq:paid-certificate}
\end{equation}
Consequently, whenever a selected grid point receives a valid token
\(m\), the whole cell \(Q_{T,m}\) is certified to lie inside the
continuous Wiener sausage. The union of the distinct valid selected
cells therefore forms a computable inner approximation of the sausage,
whose volume can be included in the exponential weight without
overestimating the true swept volume.

\subsection{Guide and candidate allocation}
In this subsection we construct the volume--confinement trade-off functions
used by the method.
The exact continuation value \(H_T^\star\) introduced in
Section~\ref{sec:ideal-law} is unavailable. We do not attempt to approximate
it pointwise. Instead, we construct a positive Lyapunov-type continuation
guide whose volume branch reflects the sequential form of the classical
Donsker--Varadhan volume--confinement balance
(Proposition~\ref{prop:target-rate}). The guide is designed to favor
the correct rare-event geometry and to satisfy the one-step inequalities used
later in the second-moment analysis. Its accuracy as an approximation to
\(H_T^\star\) is not required for exactness: the unbiasedness of the estimator
will instead follow from candidate augmentation and the empirical
proper-weight correction introduced in the next subsection.

For \(q\ge0\) and \(A\ge0\), where \(q\) represents the fraction of the time
horizon remaining and \(A\) represents the certified volume already included
in the exponential weight, set
\begin{equation}
  \JT_T(q,A)=
  \begin{cases}
    0,&q=0,\\
    \Jc_Tq^p-A,&q>0,\ 0\le A\le\widetilde v_{\star,T}q^p,\\
    q\chic_{d,T}A^{-2/d},&
       q>0,\ A>\widetilde v_{\star,T}q^p.
  \end{cases}
  \label{eq:action}
\end{equation}
For \(q>0\), the same function has the variational representation
\[
  \JT_T(q,A)
  =
  \inf_{\substack{v>0\\v\ge A}}
  \left\{
    v-A+q\chic_{d,T}v^{-2/d}
  \right\}.
\]
This is the state-dependent continuation analogue of the static
Donsker--Varadhan cost
\[
 \inf_{v>0}\left\{v+\chi_dv^{-2/d}\right\}.
\]
In the sequential expression, \(v-A\) is the additional certified volume
that would still be incurred if the eventual volume level were \(v\), while
\(q\chic_{d,T}v^{-2/d}\) is the spectral cost of confinement over the
remaining fraction \(q\) of the horizon. The constraint \(v\ge A\) reflects
that the eventual paid volume cannot be smaller than the volume already
included in the weight.

The unconstrained minimizer is
\(\widetilde v_{\star,T}q^p\). Consequently, when
\(A\le\widetilde v_{\star,T}q^p\), the optimizer remains at this preferred
volume and gives the linear branch
\(\Jc_Tq^p-A\). When \(A>\widetilde v_{\star,T}q^p\), the constraint is
active, the optimizer is \(v=A\), and only the remaining confinement cost
\(q\chic_{d,T}A^{-2/d}\) remains.

Let $\gamma_T=\frac{\Jc_T+2}{T}.$
For \(k<K_T\), define the positive volume-plus-escape guide
\begin{equation}
  H_{T,k}(x,D)
  =
  \exp\{-T\JT_T(q_k,\delta_T^d|D|)\}
  +
  \exp\{-\gamma_Td_T^{\cem}(x)\},
  \qquad
  H_{T,K_T}\equiv1.
  \label{eq:guide}
\end{equation}
Here \(D\) is a set of valid cell labels, so
\(\delta_T^d|D|\) is the certified paid volume represented by those cells.
The first summand is the volume--confinement branch derived from
\eqref{eq:action}. It supplies the rare-event guidance used to score
candidate continuations.

The second summand is a defensive escape branch. Deep inside the finite
grid window it is exponentially small, while at or outside the window
\(d_T^{\cem}(x)=0\) and the term equals one. It therefore provides a
nonvanishing guide on remote and dagger histories and is used later to
control candidates that approach or leave the finite spatial window.
The choice \(e^{-\gamma_T d_T^{\cem}(x)}\) is also convenient for the
second-moment analysis: since \(d_T^{\cem}\) is \(1\)-Lipschitz and
\(\gamma_T=O(T^{-1})\), its blockwise fluctuations accumulate only a
subexponential \(e^{o(T)}\) factor and hence do not alter the leading
second-moment rate. The terminal convention \(H_{T,K_T}\equiv1\) reflects that no continuation
guidance is needed after the final block.

We next choose the candidate-cloud size. Let
\[
  b_{T,\min}=\min_{k<K_T}b_{T,k},
  \qquad
  b_{T,\max}=\max_{k<K_T}b_{T,k},
\]
and define
\begin{equation*}
  \calL_{T,k}
  =
  T\Jc_T(q_k^p-q_{k+1}^p),
  \qquad
  C_{T,d}^{\diag}
  =
  2+4d\exp\{2\gamma_T^2b_{T,\max}\}.
\end{equation*}
The quantity \(\calL_{T,k}\) measures the maximal release of the linear
volume branch between two consecutive stages. The constant
\(C_{T,d}^{\diag}\) provides a state-uniform envelope for the squared
one-step score, including the escape branch and the Gaussian endpoint
displacement over one block.

Set
\begin{equation}
  X_{T,k}
  =
  C_{T,d}^{\diag}K_T^2e^{2\calL_{T,k}},
  \qquad
  M_{T,k}
  =
  \lceil X_{T,k}\rceil.
  \label{eq:candidate-count}
\end{equation}
Then \(M_{T,k}\) is the number of conditionally independent Brownian
candidates generated at stage \(k\). This choice is made so that the
same-candidate, or diagonal, contribution in the later second-moment
expansion is at most \(K_T^{-2}\) per block. Its accumulated effect is
therefore negligible as \(K_T\to\infty\).
Moreover, the exact finite-horizon unbiasedness does not depend on the guide
being close to \(H_T^\star\); it follows from retaining the empirical
candidate-score normalizer in the proper-weight update defined next.

\subsection{State, candidate cloud, and proper weight}

With the deterministic grids, continuation guide, and candidate counts fixed,
we now define the main candidate-augmented skeleton transition. At each stage,
the method draws a finite cloud of ordinary Brownian continuations, assigns
each candidate an immediate new-volume factor multiplied by its next-stage
guide, and selects one candidate proportionally to this score. The empirical
mean of all candidate scores is not discarded: it is retained in an
accumulated proper-weight factor. This random empirical-normalizer correction
is the augmentation mechanism that yields the exact one-step cancellation
proved in Section~\ref{sec:exactness}.

At stage \(k\), with $x_0 = 0$, write $R_k=(x_0,\ldots,x_{j_k}),$
for the selected grid path prefix. Let \(D_k\subset\calC_T\) be the set of
distinct valid positive-time tokens accumulated along this prefix, so that
\(\delta_T^d|D_k|\) is the certified paid volume at stage \(k\). Let $\mathcal G_k=(I_0,\ldots,I_{k-1})$
record the lineage of candidate labels selected in the previous blocks, and
let \(\BetaHat_{T,k}>0\) be the accumulated proper-weight factor through stage
\(k\).

The full augmented state and the reduced state used to generate and score the
next candidate cloud are
\begin{equation*}
  S_{T,k}
  =
  (k,R_k,D_k,\mathcal G_k,\BetaHat_{T,k}),
  \qquad
  S_{T,k}^{\mathrm{red}}
  =
  (k,x_{j_k},D_k).
\end{equation*}
The full state retains the selected path, label lineage, and accumulated
weight needed to define the augmented proposal law. The decision state
\(S_{T,k}^{\mathrm{red}}\) contains only the current endpoint and paid-token set needed for
the next transition. Initialize
\[
  R_0=(0),
  \qquad
  D_0=\varnothing,
  \qquad
  \mathcal G_0=(),
  \qquad
  \BetaHat_{T,0}=H_{T,0}(0,\varnothing).
\]

From the current endpoint \(y_0=x_{j_k}\), one base candidate is a block of
fine Brownian grid points $y=(y_1,\ldots,y_{n_{T,k}}).$
Its conditional product-Gaussian kernel is
\begin{equation}
  p_{T,k}(S_{T,k},\dd y)
  =
  \prod_{r=1}^{n_{T,k}}
  (2\pi h_T)^{-d/2}
  \exp\left\{
    -\frac{|y_r-y_{r-1}|^2}{2h_T}
  \right\}
  \dd y_r.
  \label{eq:base-kernel}
\end{equation}
Equivalently, with \(N(0,I_d)\) denoting the \(d\)-dimensional standard normal
distribution. 
\[
  y_r=y_{r-1}+\sqrt{h_T}Z_r,
  \qquad
  Z_r\stackrel{\mathrm{iid}}{\sim}N(0,I_d).
\]

For a candidate \(y\), define
\begin{align}
  R_k^y
    &=
    (x_0,\ldots,x_{j_k},y_1,\ldots,y_{n_{T,k}}),
    \nonumber\\
  \calE_T(y)
    &=
    \bigl\{
      \tau_T(y_r):
      1\le r\le n_{T,k},\
      \tau_T(y_r)\ne\dagger
    \bigr\},
    \nonumber\\
  \Delta D_k(y)
    &=
    \calE_T(y)\setminus D_k,
    \qquad
    D_k^y
    =
    D_k\cup\Delta D_k(y),
    \nonumber\\
  \Delta V_{T,k}(S_{T,k},y)
    &=
    \delta_T^d|\Delta D_k(y)|,
    \qquad
    x_k^y
    =
    y_{n_{T,k}}.
  \label{eq:candidate-update}
\end{align}
Here \(R_k^y\) is the grid prefix obtained if \(y\) is selected,
\(\calE_T(y)\) is the set of valid cells visited by the candidate,
\(\Delta D_k(y)\) contains only the newly certified cells, and \(D_k^y\) is
the resulting paid-token set. The quantity \(\Delta V_{T,k}(S_{T,k},y)\) is the new
certified volume contributed by the candidate, while \(x_k^y\) is its final
endpoint.

Given \(S_{T,k}\), draw \(M_{T,k}\) conditionally i.i.d. base candidates
\[
  Y_k^i
  =
  \bigl(
    Y_{k,1}^i,\ldots,Y_{k,n_{T,k}}^i
  \bigr),
  \qquad
  i=1,\ldots,M_{T,k},
\]
from the product-Gaussian kernel \eqref{eq:base-kernel}. For candidate \(i\), define
\begin{equation}
  g_{k,i}
  =
  e^{-T\Delta V_{T,k}(S_{T,k},Y_k^i)},
  \qquad
  u_{k,i}
  =
  g_{k,i}
  H_{T,k+1}
  \bigl(
    x_k^{Y_k^i},
    D_k^{Y_k^i}
  \bigr),
  \qquad
  \bar u_k
  =
  \frac{1}{M_{T,k}}
  \sum_{i=1}^{M_{T,k}}u_{k,i}.
  \label{eq:scores}
\end{equation}
The factor \(g_{k,i}\) is the candidate's immediate exponential weight for
its newly certified volume. The quantity \(u_{k,i}\) is its candidate score:
the immediate factor multiplied by the continuation guide at the resulting
state. It is the computable analogue of the ideal
``current exponential factor times future value'' score from
Section~\ref{sec:ideal-law}, but its guide component need not approximate
\(H_T^\star\) pointwise. Finally, \(\bar u_k\) is the empirical mean of the
candidate scores, and \(M_{T,k}\bar u_k\) is their categorical normalizer.

Because every candidate score is strictly positive, select a label \(I_k\)
according to
\begin{equation}
  \Psi_T
  \bigl(
    I_k=i
    \mid
    S_{T,k},Y_k^{1:M_{T,k}}
  \bigr)
  =
  \frac{u_{k,i}}{M_{T,k}\bar u_k}.
  \label{eq:selection}
\end{equation}
Thus candidates that add less certified volume and receive stronger
continuation guidance are selected more frequently. Only the selected
candidate advances the retained grid history:
\[
  R_{k+1}
  =
  R_k^{Y_k^{I_k}},
  \qquad
  D_{k+1}
  =
  D_k^{Y_k^{I_k}},
  \qquad
  \mathcal G_{k+1}
  =
  (\mathcal G_k,I_k).
\]
The unselected candidates remain part of the augmented cloud but do not alter
the selected path. In this way, an importance sampling step is done implicitly.
The accumulated proper-weight factor is updated by retaining the empirical
score mean:
\begin{equation}
  \BetaHat_{T,k+1}
  =
  \BetaHat_{T,k}
  \frac{\bar u_k}
       {H_{T,k}(x_{j_k},D_k)}.
  \label{eq:weight-update}
\end{equation}
The factor \(\bar u_k\) is the empirical-normalizer correction associated
with the score-proportional selection in \eqref{eq:selection}. Its retention
in \eqref{eq:weight-update}, rather than the accuracy of the guide, yields the
one-step proper-weight identity established later in
Lemma~\ref{lem:proper-weight}.

After the final block, the accumulated skeleton proper weight is
\begin{equation*}
  \BetaHat_T
  =
  \BetaHat_{T,K_T}
  =
  H_{T,0}(0,\varnothing)
  \prod_{k<K_T}
  \frac{\bar u_k}
       {H_{T,k}(x_{j_k},D_k)}.
\end{equation*}
Separately, the selected immediate factors multiply to the total grid
exponential payoff
\begin{equation*}
  G_T
  =
  \prod_{k<K_T}g_{k,I_k}
  =
  e^{-T\delta_T^d|D_{K_T}|}
  =
  e^{-T|\calI_T|},
  \qquad
  \calI_T
  =
  \bigcup_{m\in D_{K_T}}Q_{T,m}.
\end{equation*}
Here \(\calI_T\) is the final certified paid region: the union of all
distinct valid grid cells contributed by the selected candidate lineage.
By \eqref{eq:paid-certificate}, $\calI_T\subset W_T^{\eps_T}$
pathwise.

To recover the continuous Wiener-sausage functional, conditionally on the
selected grid prefix, complete every fine grid gap once by an ordinary
Brownian bridge. The bridges are conditionally independent across gaps and
coordinates, and the same completed path is retained for all subsequent
queries. The exact measurable construction of the
persistent completion is given in Appendix~\ref{app:bridge-init}.

\(G_T\) accounts only for the certified inner region
\(\calI_T\). To correct for the remaining region of the continuous sausage,
let \(\calN_T\) be an independent homogeneous Poisson point process on
\(\mathbb R^d\) with intensity \(T\), and define the residual void indicator
\begin{equation}
  Y_T^{\res}
  =
  \ind{
    \calN_T
    \bigl(
      W_T^{\eps_T}\setminus\calI_T
    \bigr)
    =0
  }.
  \label{eq:residual-bit}
\end{equation}
The exact finite procedure for generating this bit, using a certified bounded
path enclosure, a finite Poisson list, and persistent membership tests, is
described in Appendix~\ref{app:residual-query}.
Conditional on the completed path and the
grid data, we have
\[
  \E\!\left[
    Y_T^{\res}
    \,\middle|\,
    \text{completed path and grid data}
  \right]
  =
  \exp\!\left\{
    -T
    \bigl|
      W_T^{\eps_T}\setminus\calI_T
    \bigr|
  \right\}
  =
  \frac{F_T}{G_T}.
\]
Thus the residual bit exactly corrects the portion of the continuous sausage
that was not included in the inward grid payment. The final estimator is
\begin{equation}
    Z_T=\BetaHat_TY_T^{\res}.
  \label{eq:estimator}
\end{equation}
We summarize our importance sampling method in
Algorithm~\ref{alg:package-w}.
\begin{algorithm}[H]
\algcaption{Candidate augmentation with exact residual self-correction}
\label{alg:package-w}
\small

\textbf{Input:} fixed \(d\ge2\), \(a>0\), and \(T>0\).

\begin{enumerate}
\item Compute the deterministic constants, grids, continuation guides, and
candidate counts defined in
\eqref{eq:certified-constants}--\eqref{eq:candidate-count}.

\item Initialize $R_0=(0),$ $D_0=\varnothing,$ $\mathcal G_0=(),$ and $\BetaHat\leftarrow H_{T,0}(0,\varnothing).$

\item For \(k=0,\ldots,K_T-1\), do the following.
\begin{enumerate}
\item For each \(i=1,\ldots,M_{T,k}\), draw \(Y_k^i\) from the
product-Gaussian kernel \eqref{eq:base-kernel}. Apply the inward token map,
retain every real endpoint, and compute the candidate updates $R_k^{Y_k^i}$,
  $\Delta D_k(Y_k^i),$
  $D_k^{Y_k^i},$
  and
  $\Delta V_{T,k}(S_{T,k},Y_k^i).$

Compute \(g_{k,i}\) and \(u_{k,i}\) from \eqref{eq:scores}.

\item Compute \(\bar u_k\), and draw the selected label \(I_k\) according
to \eqref{eq:selection}.

\item Retain the selected successor state $R_{k+1}=R_k^{Y_k^{I_k}},$ $D_{k+1}=D_k^{Y_k^{I_k}},$ and $\mathcal G_{k+1}=(\mathcal G_k,I_k),$

and update the accumulated proper weight by $\BetaHat
  \leftarrow
  \BetaHat
  \frac{\bar u_k}
       {H_{T,k}(x_{j_k},D_k)}.$
\end{enumerate}

\item Conditional on the selected grid prefix, complete every grid gap
once by an ordinary Brownian bridge, retaining the same completed path for
all subsequent queries.

\item Use the finite residual-query construction in
Appendix~\ref{app:residual-query} to generate \(Y_T^{\res}\) for
\(W_T^{\eps_T}\setminus\calI_T\).
\end{enumerate}

\textbf{Output:} $Z_T=\BetaHat Y_T^{\res}.$
\end{algorithm}

We count the computational work of Algorithm~\ref{alg:package-w}
using the cost model specified in Appendix~\ref{app:work-model}; see also
\cite{BlanchetJinLiu2026Exact}.
The following proposition establishes well-definedness and shows that
the full algorithm has finite expected total work at every fixed horizon \(T\),
with at most subexponential growth in \(T\).

\begin{proposition}[Well-definedness and expected total work]
\label{prop:well-defined}
Fix \(d\ge2\) and \(a>0\).
For every \(T>0\), Algorithm~\ref{alg:package-w} defines a
probability law \(\Psi_T\) and an everywhere-defined nonnegative
measurable estimator \(Z_T\).
Let \(\mathsf{Work}_T\) denote the total work in the model of
Appendix~\ref{app:work-model}. Then
\[
  \E_{\Psi_T}\mathsf{Work}_T<\infty
  \quad\text{for every }T>0,
  \qquad
  \E_{\Psi_T}\mathsf{Work}_T=e^{o(T)}
  \quad\text{as }T\to\infty.
\]
\end{proposition}

\begin{proof}
The construction uses finitely many well-defined Borel
transition kernels with strictly positive categorical normalizers,
so the Ionescu--Tulcea theorem yields the candidate-augmented
skeleton law \(\Psi_T^{\sk}\).

Conditional on the selected grid prefix, Appendix~\ref{app:exact-real}
constructs persistent Brownian bridges, a certified box containing
the sausage, and a finite Poisson list in that box, with almost
surely finite initialization.
The list almost surely avoids the Lebesgue-null union of fine-gap
tangency sets and paid-cell faces, so every residual-membership
query terminates.
Fixed defaults on the remaining null event make \(Z_T\) everywhere
defined and measurable.
The expected-work bounds follow from
Theorem~\ref{thm:expected-work} in Appendix~\ref{app:work-bound}.
\end{proof}

\section{Exactness and finite-horizon unbiasedness}
\label{sec:exactness}

This section proves exactness at a fixed horizon. The argument has two
independent ingredients: an exact one-step cancellation on the augmented
candidate cloud and the Poisson void identity for the portion of the
continuous sausage not covered by the certified grid region.

Let \(\mathcal A_T\) denote the sigma-field generated by the candidate
clouds, selected grid prefix, completed Brownian path, and paid region
\(\calI_T\).

\begin{lemma}[Residual conditional moments]
\label{lem:residual-moments}
Conditional on \(\mathcal A_T\),
\begin{equation*}
 \E[Y_T^{\res}\mid\mathcal A_T]=\frac{F_T}{G_T},\qquad
 \E[G_TY_T^{\res}\mid\mathcal A_T]=F_T,\qquad
 \E[(G_TY_T^{\res})^2\mid\mathcal A_T]=F_TG_T.
\end{equation*}
\end{lemma}

\begin{proof}
\eqref{eq:paid-certificate} implies
\(\calI_T\subset W_T^{\eps_T}\) pathwise. By the definition
\eqref{eq:residual-bit} and independence of \(\calN_T\), the conditional
Poisson void probability is
\[
 \E[Y_T^{\res}\mid\mathcal A_T]
 =\exp\{-T|W_T^{\eps_T}\setminus\calI_T|\}
 =\exp\{-T|W_T^{\eps_T}|+T|\calI_T|\}
 =\frac{F_T}{G_T}.
\]
The remaining identities follow from \(G_T=e^{-T|\calI_T|}\) and
\((Y_T^{\res})^2=Y_T^{\res}\).
\end{proof}

Recall that \(\Psi_T^{\sk}\) is the marginal proposal law of the
candidate clouds, selected labels, selected grid prefix, and
accumulated proper weight, before the Brownian-bridge and
residual-Poisson extensions. We next establish the one-step weighting identity used in the
unbiasedness proof.

\begin{lemma}[One-step proper weighting]
\label{lem:proper-weight}
Fix a stage \(k\), a reachable full state $S=(k,R,D,\mathcal G,\beta),$
and let \(x\) be the terminal point of \(R\). For a nonnegative measurable
function \(f\) of the selected candidate,
\begin{equation*}
\begin{split}
 \E_{\Psi_T}\left[
   \left.
   \frac{\bar u_k}{H_{T,k}(x,D)}f(Y_k^{I_k})
   \right|S_{T,k}=S
 \right] =
 \int
 \frac{
   e^{-T\Delta V_{T,k}(S,y)}
   H_{T,k+1}(x_k^y,D_k^y)
 }{H_{T,k}(x,D)}
 f(y)\,p_{T,k}(S,\dd y).
\end{split}
\end{equation*}
\end{lemma}

\begin{proof}
Write \(M=M_{T,k}\). Conditional on the iid candidate cloud and on
\(S_{T,k}=S\), the selection rule \eqref{eq:selection} gives
\[
 \E\left[
   \left.
   \frac{\bar u_k}{H_{T,k}(x,D)}f(Y_k^{I_k})
   \right|Y_k^{1:M},S_{T,k}=S
 \right]
 =\frac1M\sum_{i=1}^M
 \frac{u_{k,i}}{H_{T,k}(x,D)}f(Y_k^i).
\]
The empirical normalizer cancels exactly. Averaging the iid summands and
using the definition of \(u_{k,i}\) in \eqref{eq:scores} finishes the proof.
\end{proof}

\begin{corollary}[Global proper weighting]
\label{cor:global-proper}
For every nonnegative measurable functional \(\varphi\) of the selected
grid prefix,
\begin{equation*}
 \E_{\Psi_T^{\sk}}[\BetaHat_T\varphi(R_{K_T})]
 =\E_{\overline\Pp_T^{\grid}}[G_T\varphi(R_{K_T})],
\end{equation*}
where \(\overline\Pp_T^{\grid}\) is the ordinary Brownian endpoint law on
the \(N_T\)-gap grid.
\end{corollary}

\begin{proof}
The candidate kernel and scores depend on the past only through the selected
grid history. Starting from
\(\BetaHat_{T,0}=H_{T,0}(0,\varnothing)\), backward induction using
Lemma~\ref{lem:proper-weight} therefore integrates out one cloud and its
selected label at a time. Along an ordinary Brownian grid path write $g_k=\exp\{-T\delta_T^d(|D_{k+1}|-|D_k|)\}.$
These are the same incremental factors as \(g_{k,I_k}\) along a selected
proposal path. The induction gives
\[
\begin{aligned}
 \E_{\Psi_T^{\sk}}[\BetaHat_T\varphi(R_{K_T})]
 &=\E_{\overline\Pp_T^{\grid}}\!\left[
    H_{T,0}(0,\varnothing)
    \prod_{k<K_T}
    \frac{g_kH_{T,k+1}(x_{j_{k+1}},D_{k+1})}
         {H_{T,k}(x_{j_k},D_k)}
    \varphi(R_{K_T})
   \right]\\
 &=\E_{\overline\Pp_T^{\grid}}[G_T\varphi(R_{K_T})].
\end{aligned}
\]
Indeed, the intermediate guide factors telescope,
\(H_{T,K_T}\equiv1\), and \(\prod_{k<K_T}g_k=G_T\). The balanced blocks
partition the \(N_T\) Gaussian grid gaps, so the iterated base kernels
generate exactly \(\overline\Pp_T^{\grid}\).
\end{proof}

The same conditional Brownian-bridge kernel is used to complete a selected
proposal grid path or an ordinary Brownian grid path. Combining this common
completion with Corollary~\ref{cor:global-proper} and the residual void
identity proves unbiasedness.

\begin{theorem}[Exact finite-horizon unbiasedness]
\label{thm:unbiasedness}
For every fixed \(d\ge2\), \(a>0\), and \(T>0\), the estimator in
\eqref{eq:estimator} is nonnegative, integrable, and
\begin{equation*}
 \E_{\Psi_T}Z_T
 =\E_{\Pp_T}\exp\{-T|W_T^{aT^{-1/d}}|\}
 =\alpha_T.
\end{equation*}
\end{theorem}

\begin{proof}
By Proposition~\ref{prop:well-defined}, \(Z_T\) is a nonnegative measurable
random variable under the well-defined law \(\Psi_T\). All expectations
below may therefore be evaluated by Tonelli's theorem before integrability
is known.

For a grid prefix \(R\), let $\varphi(R)
  =\E_{\mathrm{br}}\!\left[
      \left.\frac{F_T}{G_T}\,\right|R
    \right],$
where \(\E_{\mathrm{br}}\) integrates the ordinary conditional
Brownian-bridge completion. The paid-set certificate implies
\(0<\varphi\le1\). Then, by Lemma~\ref{lem:residual-moments} and  Corollary~\ref{cor:global-proper}, we have
\[
\begin{aligned}
  \E_{\Psi_T}Z_T
  &=\E_{\Psi_T^{\sk}}[\BetaHat_T\varphi(R_{K_T})]=\E_{\overline\Pp_T^{\grid}}[G_T\varphi(R_{K_T})]
   =\E_{\Pp_T}F_T
   =\alpha_T.
\end{aligned}
\]
\end{proof}

\section{Second-moment decomposition}
\label{sec:second-moment}

Exact unbiasedness determines the first moment of the estimator. To assess
its rare-event efficiency, we compare the second moment with the smallest
exponential scale permitted by unbiasedness.
We use the standard asymptotic optimality criterion of rare-event simulation
\cite{Siegmund1976SequentialTests,GlynnWhitt1992Efficiency}.

\begin{definition}[Asymptotic optimality]
\label{def:log-eff}
An unbiased estimator \(Z_T\) of \(\alpha_T\), generated under
\(\Psi_T\), is \emph{asymptotically optimal} at speed \(T\) (also called
logarithmically efficient) if
\[
 \limsup_{T\to\infty}\frac1T
 \log\frac{\E_{\Psi_T}Z_T^2}{\alpha_T^2}\le0.
\]
\end{definition}

By Proposition~\ref{prop:target-rate},
\(\log\alpha_T=-JT+o(T)\). Moreover, exact unbiasedness and Jensen's
inequality imply
\[
 \E_{\Psi_T}Z_T^2\ge\alpha_T^2=e^{-2JT+o(T)}.
\]
Consequently, asymptotic optimality is equivalent in the present setting to
the upper bound
\[
 \log\E_{\Psi_T}Z_T^2\le-2JT+o(T).
\]
We now expand the second moment itself. This expansion will identify the
one-step operator and the diagonal and off-diagonal quantities that must be
controlled.

\subsection{From the second moment to the one-step operator}

For the second-moment recursion, retain the selected grid history $\Xi_k=(k,R_k,D_k).$
We do not retain the selected-label history
\(\mathcal G_k=(I_0,\ldots,I_{k-1})\) or the accumulated proper weight
\(\widehat\beta_{T,k}\) as components of \(\Xi_k\). The randomness of the
past labels is already encoded in the realized selected grid prefix \(R_k\)
and paid set \(D_k\), while the averaging over the current selected label is
performed explicitly in the one-step conditional expectation below. The
multiplicative evolution of \(\widehat\beta_{T,k}\) is peeled off one stage at
a time. On the other hand, the full prefix \(R_k\) must be retained,
rather than only the decision state \((k,x_{j_k},D_k)\), because the terminal
residual factor after Brownian-bridge completion depends on the entire
selected grid history.

For \(\xi=(k,R,D)\), write \(x(\xi)\) for the terminal point of \(R\), and
use the shorthand $H_{T,k}(\xi):=H_{T,k}(x(\xi),D).$
Define the terminal bridge-completion factor
\begin{equation*}
 \phi_T(\Xi_{K_T})
 :=\E\left[
   \left.\frac{F_T}{G_T}\right|R_{K_T}
 \right],
 \qquad 0\le\phi_T\le1,
\end{equation*}
where the conditional expectation is under the ordinary Brownian-bridge
completion of the selected grid prefix. The bound follows from
\(\calI_T\subset W_T^{\eps_T}\), which implies \(0<F_T/G_T\le1\).

Set $C_{T,k}:=\frac{\bar u_k}{H_{T,k}(\Xi_k)}.$
By the proper-weight recursion,
\[
 \BetaHat_T
 =H_{T,0}(0,\varnothing)\prod_{k<K_T}C_{T,k}.
\]
Let \(\mathcal A_T^{\sk}\) be the sigma-field generated by all candidate
clouds and selected labels. The proper weight and selected grid prefix are
\(\mathcal A_T^{\sk}\)-measurable. Conditional on that prefix, the ordinary
bridge completion is independent of the remaining augmentation, and the
residual Poisson process is independent of both. Since \(Y_T^{\res}\) is an
indicator, Lemma~\ref{lem:residual-moments} therefore gives
\[
 \E\bigl[(Y_T^{\res})^2\mid\mathcal A_T^{\sk}\bigr]
 =\E[Y_T^{\res}\mid\mathcal A_T^{\sk}]
 =\phi_T(\Xi_{K_T}).
\]
Multiplying by \(\BetaHat_T^2\) and integrating yields the exact skeleton
representation
\begin{equation}
 \E_{\Psi_T}Z_T^2
 =H_{T,0}(0,\varnothing)^2
 \E_{\Psi_T^{\sk}}\left[
   \left(\prod_{k<K_T}C_{T,k}^2\right)
   \phi_T(\Xi_{K_T})
 \right].
 \label{eq:raw-second-moment}
\end{equation}

For one fresh base candidate \(Y\), let
\(\xi^Y=(k+1,R^Y,D^Y)\) be the successor history obtained by appending \(Y\)
to \(R\) and updating the paid-token set as in
\eqref{eq:candidate-update}.
Equation~\eqref{eq:raw-second-moment} shows what must be propagated backward.
For every nonnegative measurable function \(f\) on stage-\(k+1\) histories,
define the one-step squared-weight operator
\begin{equation}
 (\calK_{T,k}f)(\xi)
 :=
 \E_{\Psi_T^{\sk}}\!\left[
   \left(\frac{\bar u_k}{H_{T,k}(\xi)}\right)^2
   f(\Xi_{k+1})
   \,\middle|\,
   \Xi_k=\xi
 \right].
 \label{eq:pair-operator-definition}
\end{equation}
Under this one-step proposal transition,
\(\Xi_{k+1}=\xi^{Y_k^{I_k}}\). Applying the tower property first at stage
\(K_T-1\) and then backward through the remaining stages gives
\begin{equation}
 \E_{\Psi_T}Z_T^2
 =H_{T,0}(0,\varnothing)^2
 \bigl(\calK_{T,0}\calK_{T,1}\cdots
       \calK_{T,K_T-1}\phi_T\bigr)(\Xi_0).
 \label{eq:exact-second-operator}
\end{equation}

\subsection{Same-cloud diagonal--off-diagonal decomposition}
In this subsection we further decompose \(\calK_{T,k}\). Since the candidate volume update depends on
the full state only through \(\xi\), we write \(\Delta V_{T,k}(\xi,Y)\) for the same
new paid volume. Define
\begin{align*}
 g_{T,k}(\xi,Y)
 &:=\exp\{-T\Delta V_{T,k}(\xi,Y)\},\nonumber\\
 W_{T,k}(\xi,Y)
 &:=\frac{g_{T,k}(\xi,Y)H_{T,k+1}(\xi^Y)}{H_{T,k}(\xi)}.
\end{align*}
Here and below, expectations over fresh candidates are under the base kernel
\eqref{eq:base-kernel}, conditional on the current history \(\xi\).

\begin{lemma}[Same-cloud ordered-pair representation]
\label{lem:pair-operator}
For every nonnegative measurable \(f\), the operator in
\eqref{eq:pair-operator-definition} satisfies
\begin{equation}
\begin{split}
 (\calK_{T,k}f)(\xi)
 &=
 \frac1{M_{T,k}}\E[W_{T,k}(\xi,Y)^2f(\xi^Y)]+
 \left(1-\frac1{M_{T,k}}\right)
 \E[W_{T,k}(\xi,Y)W_{T,k}(\xi,Y')f(\xi^Y)].
\end{split}
\label{eq:pair-operator}
\end{equation}
where \(Y\) and \(Y'\) are independent base candidates from the same current
history, and only \(Y\) determines the successor primary history. Define
\begin{equation*}
 \mathfrak{d}_{T,k}(\xi):=\E[W_{T,k}(\xi,Y)^2],
 \qquad
 \rho_{T,k}(\xi):=\E[W_{T,k}(\xi,Y)].
\end{equation*}
Then the total mass of the operator is
\begin{equation}
 m_{T,k}(\xi)
 :=(\calK_{T,k}\mathbf 1)(\xi)
 =\frac{\mathfrak{d}_{T,k}(\xi)}{M_{T,k}}
 +\left(1-\frac1{M_{T,k}}\right)\rho_{T,k}(\xi)^2.
 \label{eq:pair-mass}
\end{equation}
\end{lemma}

\begin{proof}
Abbreviate \(M=M_{T,k}\), \(H=H_{T,k}(\xi)\), and
\(u_i=u(Y^i)\). Conditional on the candidate cloud
\(Y^{1:M}=(Y^1,\ldots,Y^M)\), the score-proportional selection rule gives
\[
 \Psi_T^{\sk}(I=i\mid Y^{1:M},\xi)=\frac{u_i}{M\bar u}.
\]
Consequently,
\begin{align*}
 \E_{\Psi_T^{\sk}}\!\left[
   \left(\frac{\bar u}{H}\right)^2f(\xi^{Y^I})
   \,\middle|\,Y^{1:M},\xi
 \right]=
 \sum_{i=1}^M\frac{u_i}{M\bar u}
 \left(\frac{\bar u}{H}\right)^2 f(\xi^{Y^i})
 =\frac{\bar u}{MH^2}\sum_{i=1}^M u_i f(\xi^{Y^i}).
\end{align*}
Using \(\bar u=M^{-1}\sum_{j=1}^M u_j\) and
\(W_{T,k}(\xi,Y^i)=u_i/H\), the last expression becomes
\[
 \frac1{M^2}\sum_{i,j=1}^M
 W_{T,k}(\xi,Y^i)W_{T,k}(\xi,Y^j)f(\xi^{Y^i}).
\]
There are \(M\) diagonal terms \(i=j\), with total coefficient \(1/M\),
and \(M(M-1)\) off-diagonal ordered terms \(i\ne j\), with total
coefficient \(1-1/M\). Averaging over the iid cloud and using independence
of distinct cloud coordinates in the off-diagonal terms gives
\eqref{eq:pair-operator}. Taking \(f\equiv1\), and using independence to
write
\(\E[W_{T,k}(\xi,Y)W_{T,k}(\xi,Y')]=\rho_{T,k}(\xi)^2\), gives
\eqref{eq:pair-mass}.
\end{proof}

The decomposition \eqref{eq:pair-mass} exposes the two distinct parts of the
second-moment problem. The same-candidate, or diagonal, contribution is
\(\mathfrak{d}_{T,k}/M_{T,k}\); it is controlled by the candidate multiplicity. The
different-candidate, or off-diagonal, contribution is governed by
\(\rho_{T,k}^2\); its control requires the blockwise lower-volume and
continuation analysis in Section~\ref{sec:efficiency}.

\subsection{Diagonal control and normalized global representation}

Write
\(H_{T,k}=H_{T,k}^{\mathrm{vol}}+H_{T,k}^{\esc}\) for the two summands in
\eqref{eq:guide}. At the terminal stage set
\(H_{T,K_T}^{\mathrm{vol}}\equiv1\) and
\(H_{T,K_T}^{\esc}\equiv0\), so the same decomposition remains valid when
\(k=K_T-1\). Accordingly, decompose
\[
 W_{T,k}(\xi,Y)
 =W_{T,k}^{\mathrm{vol}}(\xi,Y)
 +W_{T,k}^{\esc}(\xi,Y),
\]
where
\[
 W_{T,k}^{\mathrm{vol}}(\xi,Y)
 :=\frac{g_{T,k}(\xi,Y)H_{T,k+1}^{\mathrm{vol}}(\xi^Y)}
         {H_{T,k}(\xi)},
 \qquad
 W_{T,k}^{\esc}(\xi,Y)
 :=\frac{g_{T,k}(\xi,Y)H_{T,k+1}^{\esc}(\xi^Y)}
         {H_{T,k}(\xi)}.
\]

\begin{lemma}[State-uniform diagonal envelope]
\label{lem:diagonal-envelope}
For every stage and every reachable history,
\begin{equation*}
 \mathfrak{d}_{T,k}(\Xi_k)
 \le C_{T,d}^{\diag}e^{2\calL_{T,k}},\qquad
 \frac{\mathfrak{d}_{T,k}(\Xi_k)}{M_{T,k}}\le K_T^{-2}.
\end{equation*}
\end{lemma}

\begin{proof}
Fix a current history \(\xi=(k,R,D)\), a candidate \(Y\), and put
\[
 A=\delta_T^d|D|,
 \qquad
 w=\Delta V_{T,k}(\xi,Y),
 \qquad
 V=A+w,
 \qquad
 q=q_k,
 \qquad
 q'=q_{k+1}.
\]
Define
\[
 \Phi_T(s,B)
 :=
 B+\JT_T(s,B)
 =
 \inf_{\substack{v>0\\v\ge B}}\{v+s\chic_{d,T}v^{-2/d}\}.
\]
Since \(V\ge A\), \(\Phi_T(q,A)\le\Phi_T(q,V)\). For
\(s>0\), the minimizer in \(\Phi_T(s,V)\) is $v_s=\max\{V,\widetilde v_{\star,T}s^p\}.$
The coupled identities in \eqref{eq:certified-constants} therefore give, for
almost every \(s>0\),
\[
 \frac{\partial}{\partial s}\Phi_T(s,V)
 =
 \chic_{d,T}v_s^{-2/d}
 \le
 \chic_{d,T}(\widetilde v_{\star,T}s^p)^{-2/d}
 =
 p\Jc_Ts^{p-1}.
\]
Integrating from \(q'\) to \(q\), with the case \(q'=0\) obtained by
continuity, yields
\[
 \begin{aligned}
 \JT_T(q,A)-w-\JT_T(q',V)
 &=
 \Phi_T(q,A)-\Phi_T(q',V)\le
 \Jc_T\bigl(q^p-(q')^p\bigr)
 =
 \calL_{T,k}/T.
 \end{aligned}
\]
Exponentiating gives
\[
 e^{-Tw}H_{T,k+1}^{\mathrm{vol}}(\xi^Y)
 \le e^{\calL_{T,k}}H_{T,k}^{\mathrm{vol}}(\xi).
\]
Since the full denominator \(H_{T,k}(\xi)\) dominates its volume summand,
the normalized volume contribution to \(W_{T,k}(\xi,Y)\) is at most
\(e^{\calL_{T,k}}\).

Let \(Z\sim N(0,b_{T,k}I_d)\) be the candidate's total endpoint
displacement. The cemetery distance is one-Lipschitz in the
\(\ell^\infty\) norm. Since \(g_{T,k}\le1\) and the full denominator also
dominates its escape summand,
\[
 \frac{g_{T,k}(\xi,Y)H_{T,k+1}^{\esc}(x(\xi)+Z)}
      {H_{T,k}(\xi)}
 \le e^{\gamma_T\|Z\|_\infty}.
\]
A union bound over coordinates and Gaussian square completion gives $ \E e^{2\gamma_T\|Z\|_\infty}
 \le2d\,e^{2\gamma_T^2b_{T,k}}.$
Applying \((r+s)^2\le2r^2+2s^2\), using
\(b_{T,k}\le b_{T,\max}\), and recalling the definition of
\(C_{T,d}^{\diag}\) proves the first inequality. The second follows from
\(M_{T,k}\ge C_{T,d}^{\diag}K_T^2e^{2\calL_{T,k}}\).
\end{proof}

Thus the candidate multiplicity suppresses the same-candidate contribution
to at most \(K_T^{-2}\) per block. Lemma~\ref{lem:diagonal-envelope} also
shows that \(\mathfrak{d}_{T,k}(\xi)<\infty\). The Cauchy--Schwarz inequality gives
\(\rho_{T,k}(\xi)^2\le \mathfrak{d}_{T,k}(\xi)\), and positivity of every score gives $0<m_{T,k}(\xi)<\infty.$
The weighted one-step law in
\eqref{eq:pair-operator-definition} therefore defines a finite kernel
\(\calK_{T,k}(\xi,\dd\xi')\) on successor histories. Normalize it by
\[
 \Pi_{T,k}(\xi,\dd\xi')
 :=\frac{\calK_{T,k}(\xi,\dd\xi')}{m_{T,k}(\xi)},
\]
and let \(\Pi_T\) be the iterated one-primary history law started from
\(\Xi_0=(0,(0),\varnothing)\). 

\begin{theorem}[Exact finite-\(T\) second moment]
\label{thm:exact-second}
The operator representation \eqref{eq:exact-second-operator} admits the
normalized form
\begin{equation}
 \E_{\Psi_T}Z_T^2
 =H_{T,0}(0,\varnothing)^2
 \E_{\Pi_T}\left[
   \phi_T(\Xi_{K_T})\prod_{k<K_T}m_{T,k}(\Xi_k)
 \right].
 \label{eq:exact-second-moment}
\end{equation}
\end{theorem}

\begin{proof}
For every nonnegative measurable function \(f\),
\[
 (\calK_{T,k}f)(\xi)
 =
 m_{T,k}(\xi)
 \int f(\xi')\,\Pi_{T,k}(\xi,\dd\xi').
\]
Applying this first to the leftmost operator in \eqref{eq:exact-second-operator}
gives
\[
\begin{aligned}
(\calK_{T,0}\cdots\calK_{T,K_T-1}\phi_T)(\Xi_0)=
m_{T,0}(\Xi_0)
\int
(\calK_{T,1}\cdots\calK_{T,K_T-1}\phi_T)(\xi_1)
\,\Pi_{T,0}(\Xi_0,\dd\xi_1).
\end{aligned}
\]
Repeating the same substitution through stages
\(1,\ldots,K_T-1\) extracts the successive factors
\(m_{T,k}(\Xi_k)\) and evolves the primary history according to
\(\Pi_{T,k}\). Hence
\[
 (\calK_{T,0}\cdots\calK_{T,K_T-1}\phi_T)(\Xi_0)
 =
 \E_{\Pi_T}\left[
   \phi_T(\Xi_{K_T})
   \prod_{k<K_T}m_{T,k}(\Xi_k)
 \right],
\]
which together with \eqref{eq:exact-second-operator} proves the claim.
\end{proof}

The exact identity now exposes the full proof target. The terminal factor is
harmless because \(0\le\phi_T\le1\), and the diagonal term has been reduced
to \(K_T^{-2}\) per block by Lemma~\ref{lem:diagonal-envelope}. It remains
to identify the exponent of the initial guide and to prove the state-uniform
off-diagonal bound
\[
 \rho_{T,k}(\Xi_k)\le e^{\eta_{T,k}},
 \qquad
 \sum_{k<K_T}\eta_{T,k}=o(T).
\]

\section{Asymptotic optimality}
\label{sec:efficiency}

Our target is \(\log\E_{\Psi_T}Z_T^2\le-2JT+o(T)\). By the exact
representation in Section~\ref{sec:second-moment}, it remains to show that
the initial guide has logarithm \(-JT+o(T)\) and that the logarithmic
errors in the off-diagonal continuation means sum to \(o(T)\). We first
identify the initial scale, then translate a fixed-radius sausage
lower-tail estimate into a uniform fresh-block grid-volume bound. The
scalar continuation inequality combines this bound with four volume
sectors; escape and terminal-block estimates complete the error ledger.
We then insert these estimates into the second-moment representation.

\subsection{Initial guide}

\begin{lemma}[Initial guide]
\label{lem:initial-guide}$ \log H_{T,0}(0,\varnothing)=-JT+o(T).$
\end{lemma}

\begin{proof}
At \(q_0=1\) and \(D_0=\varnothing\), the volume branch equals
\(e^{-T\Jc_T}\). By \eqref{eq:app-cem-origin}, the inward-window distance
from the origin is \(T^2+o(T^2)\), and hence
\[
 \gamma_Td_T^{\cem}(0)=(\Jc_T+2)T+o(T).
\]
The escape branch is therefore exponentially smaller than the volume branch.
The certification error in \eqref{eq:certified-constants} gives
\(\Jc_T=J+o(1)\), proving the claim.
\end{proof}

\subsection{Fixed-radius lower-volume input and block translation}

The downward large deviations of the Wiener-sausage volume were developed by
Donsker and Varadhan, Bolthausen, and Sznitman
\cite{DonskerVaradhan1975WienerSausage,Bolthausen1990Volume,
Sznitman1990LifschitzI,Sznitman1990LifschitzII}. We use the convenient
formulation in equations (1.4)--(1.5) of van den Berg, Bolthausen, and den
Hollander \cite{VanDenBergBolthausenDenHollander2001Moderate}. For a fixed
radius \(r>0\), Brownian generator \(\Delta/2\), and a deterministic function
\(f(s)\to\infty\), it states that
\begin{equation}
 \frac{f(s)^{2/d}}{s}
 \log\Pp\{|W_s^{r,\mathrm{open}}|\le f(s)\}
 \longrightarrow-\frac{\lambda_d}{2},
 \label{eq:hard-lower-tail}
\end{equation}
provided \(f(s)=o(s/\log s)\) in \(d=2\) and \(f(s)=o(s)\) in
\(d\ge3\). Here
\[
  W_s^{r,\mathrm{open}}
  :=
  \bigcup_{0\le u\le s}B_2(B_u,r)
\]
is the open-ball sausage. Since
\(W_s^{r,\mathrm{open}}\subset W_s^r\), where \(W_s^r\) is the closed-ball
sausage used throughout this paper,
\[
  \{|W_s^r|\le f(s)\}
  \subset
  \{|W_s^{r,\mathrm{open}}|\le f(s)\}.
\]
Hence the upper bound furnished by \eqref{eq:hard-lower-tail} also applies to
the closed-ball lower-tail probabilities used below; this is the only
direction needed in the proof.

Set \(\widetilde\lambda_T(v):=\chic_{d,T}v^{-2/d}\) for \(v>0\), with
\(\widetilde\lambda_T(0):=+\infty\).
The block schedule satisfies, for fixed \(d,a\),
\begin{equation}
 b_T=o(T),\quad K_T\to\infty,\quad K_T\log N_T=o(T),\quad
 b_{T,\min}\asymp b_{T,\max}\asymp b_T.
 \label{eq:schedule-orders}
\end{equation}
Moreover,
\(b_T\delta_T^{d-2}\to\infty\) for \(d\ge3\), whereas $\frac{b_T}{\log(b_{T,\max}\delta_T^{-2})}\to\infty$ when $d=2$.
These consequences of \eqref{eq:raster-schedule}--\eqref{eq:block-rate} are
verified in Appendix~\ref{app:lower-volume}; they place every balanced block
within the source range of \eqref{eq:hard-lower-tail}.

For a fresh candidate define the infinite-grid volume, before inward-window
replacement and before previously paid tokens are removed, by
\begin{equation*}
 U_{T,k}^{\raw}
 =\delta_T^d\left|\left\{
 \lfloor Y_r/\delta_T\rfloor:1\le r\le n_{T,k}
 \right\}\right|.
\end{equation*}

The next proposition is the uniform blockwise form of this fixed-radius input:
a freeze-then-wait diagonal argument converts the pointwise source theorem
into a history- and phase-uniform lower-tail bound for the raw grid volume
defined above.

\begin{proposition}[Uniform fresh-block lower-volume bound]
\label{prop:fresh-tail}
There are deterministic step functions
\[
 u_T,\theta_T,\kappa_T,\zeta_T\downarrow0,\qquad
 A_T^{\cut},v_T^{\up}\uparrow\infty,\qquad
 v_T^{\up}=A_T^{\cut}+J+3,
\]
and nonnegative remainders \(r_{T,k}\) such that, uniformly over bulk blocks
\(q_{k+1}\ge\theta_T\), incoming histories, grid phases, and
\(v\in[u_T,v_T^{\up}]\),
\begin{equation*}
 \log\Pp\{U_{T,k}^{\raw}\le v\mid\Xi_k\}
 \le-b_{T,k}\widetilde\lambda_T(v)+r_{T,k}.
\end{equation*}
The choices may additionally be made so that
\begin{equation}
 u_T\le\frac14v_\star\theta_T^p,
 \qquad
 \sup_{k<K_T}\frac{r_{T,k}}{b_T}\longrightarrow0,
 \qquad
 \sum_{k<K_T}r_{T,k}=o(T).
 \label{eq:fresh-tail-coupling}
\end{equation}
The choices can be made so that
\begin{equation}
 \begin{split}
 K_T\log(N_T+3)+T\widetilde\lambda_T(A_T^{\cut})
 +T\theta_T^p+TK_T^{-p}+K_T+T\gamma_T^2+T|\Jc_T-J|=o(T).
 \end{split}
 \label{eq:main-ledger}
\end{equation}
\end{proposition}

\begin{proof}
The proof consists of three reductions. First, for analysis only, complete
the Gaussian endpoint candidate with independent ordinary conditional
bridges. The completed block is an ordinary Brownian increment. Integrating
out these bridges leaves the candidate kernel unchanged, while translation
removes the incoming endpoint; only the phase of the \(\delta_T\)-grid
remains.

Second, apply Lemma~\ref{lem:raster-source-translation}. For every fixed
\(\kappa,\zeta>0\), that lemma constructs a good mesh event
\(\Omega_{T,k}^{\mathrm{mesh}}(\kappa)\) and gives, uniformly over incoming histories and grid
phases,
\begin{equation}
\begin{split}
 \Pp\{U_{T,k}^{\raw}\le v\mid\Xi_k\}
 &\le
 4dN_T
 \exp\left\{
   -\frac{\kappa^2\delta_T^2}{C_dh_T}
 \right\}+
 \Pp\left\{
   |W_{s_{T,k}}^\zeta|
   \le f_{T,k}(v)
 \right\}.
\end{split}
\label{eq:main-raw-to-source}
\end{equation}
where $a_{\kappa,\zeta}
 =
 [1+2(\kappa+\zeta)]^{-d},$ $s_{T,k}
 =
 b_{T,k}\delta_T^{-2},$ and $f_{T,k}(v)
 =
 \frac{v}{a_{\kappa,\zeta}\delta_T^d}.$
The same lemma identifies the fixed-radius source exponent as
\begin{equation}
 \frac{\lambda_d}{2}
 s_{T,k}f_{T,k}(v)^{-2/d}
 =
 b_{T,k}a_{\kappa,\zeta}^{2/d}\Lambda_{\FK}(v).
 \label{eq:main-source-exponent}
\end{equation}

For this step, keep \(\kappa,\zeta>0\) and a compact threshold interval
\([u,R]\subset(0,\infty)\) fixed. The schedule relations
\eqref{eq:schedule-orders} place \((s_{T,k},f_{T,k}(v))\) in the
appropriate dimension-dependent source range of
\eqref{eq:hard-lower-tail}. The source theorem and
\eqref{eq:main-source-exponent} give an exponential upper bound
\[
 \Pp\!\left\{|W_{s_{T,k}}^\zeta|\le f_{T,k}(v)\right\}
 \le B_{T,k}(v)
 :=\exp\!\left\{
 -b_{T,k}a_{\kappa,\zeta}^{2/d}\Lambda_{\FK}(v)+o(b_T)
 \right\}.
\]
for the source probability, uniformly over balanced blocks and
\(v\in[u,R]\); the frozen-parameter uniformity is justified in
Lemma~\ref{lem:freeze-wait-uniformization}. Since \(\delta_T^2/h_T\asymp T(\log T)^3\) and \(b_T=o(T)\),
the first term in \eqref{eq:main-raw-to-source} is negligible
relative to \(B_{T,k}(v)\), uniformly over balanced blocks and
\(v\in[u,R]\). Hence
\[
 \log\Pp\{U_{T,k}^{\raw}\le v\mid\Xi_k\}
 \le-b_{T,k}a_{\kappa,\zeta}^{2/d}\Lambda_{\FK}(v)+o(b_T).
\]
For the \(T\)-dependent parameters and threshold interval in the
proposition, the waiting condition
\eqref{eq:app-mesh-domination} ensures the same domination.

Since \(a_{\kappa,\zeta}^{2/d}\to1\) as \(\kappa,\zeta\downarrow0\) and
\(\widetilde\lambda_T(v)\to\Lambda_{\FK}(v)\) on every fixed compact
interval bounded away from zero, the remaining tasks are to control the
dilation and spectral-certification losses uniformly and to upgrade the
fixed-parameter lower-tail estimate above to the \(T\)-dependent ranges
required in the proposition.

Third, apply Lemma~\ref{lem:freeze-wait-uniformization}, whose inputs are
the grid-to-source translation in
Lemma~\ref{lem:raster-source-translation}, the classical fixed-radius
lower-volume theorem \eqref{eq:hard-lower-tail}, and the schedule relations
\eqref{eq:schedule-orders}. The lemma chooses the auxiliary parameters
sufficiently slowly to obtain a bound that is uniform over bulk blocks,
incoming histories, grid phases, and \(v\in[u_T,v_T^{\up}]\), while
absorbing the dilation, spectral-certification, and uniformization losses
into nonnegative block remainders satisfying $\sup_{k<K_T}\frac{r_{T,k}}{b_T}\longrightarrow0$ and $\sum_{k<K_T}r_{T,k}=o(T).$
The same construction produces the deterministic step functions appearing
in the proposition and simultaneously enforces the error ledger
\eqref{eq:main-ledger}. This proves the proposition.
\end{proof}

\subsection{Off-diagonal continuation control}
Let \(q=q_k\), \(q'=q_{k+1}\),
\(A=\delta_T^d|D_k|\), let \(w\) be the new paid volume, and set
\(V=A+w\).
The next deterministic inequality links the one-block lower-tail cost to the
continuation guide: the immediate new-volume cost, the spectral cost over the
current block, and the continuation cost after the block together dominate
the guide value at the incoming state.

\begin{lemma}[Scalar continuation inequality]
\label{lem:scalar-continuation}
For every \(A,w\ge0\) and \(q>q'\ge0\),
\begin{equation*}
 w+(q-q')\widetilde\lambda_T(V)+\JT_T(q',V)
 \ge\JT_T(q,A).
\end{equation*}
\end{lemma}

\begin{proof}
If \(V=0\), the inequality follows from
\(\widetilde\lambda_T(0)=+\infty\). Suppose \(V>0\). For every
\(v\ge V\), the identity \(w+v-V=v-A\) and monotonicity of
\(\widetilde\lambda_T\) give
\[
 \begin{aligned}
w+(q-q')\widetilde\lambda_T(V)+(v-V)
 +q'\widetilde\lambda_T(v)\ge v-A+q\widetilde\lambda_T(v).
 \end{aligned}
\]
Taking the infimum on the left over \(v\ge V\) yields
\[
 w+(q-q')\widetilde\lambda_T(V)+\JT_T(q',V).
\]
On the right, replacing the constraint \(v\ge V\) by the weaker constraint
\(v\ge A\) can only decrease the infimum, thus the proof is finished.
\end{proof}

The preceding lower-tail estimate and scalar continuation inequality now
combine to control the off-diagonal one-step factor uniformly, with only a
subexponential error accumulated over all blocks.

Set
\[
 \mathcal B_T:=\{k<K_T:q_{k+1}\ge\theta_T\},
 \qquad
 \mathcal T_T:=\{k<K_T:q_{k+1}<\theta_T\}.
\]
For a sufficiently large dimension-dependent constant \(C_d\), define
\begin{equation}
 \eta_{T,k}^{\mathrm{bulk}}
 :=r_{T,k}+\log(N_T+3)
 +b_{T,k}\widetilde\lambda_T(A_T^{\cut})
 +C_d\gamma_T^2b_{T,k}+C_d,
 \qquad k\in\mathcal B_T,
 \label{eq:bulk-eta}
\end{equation}
and
\begin{equation}
 \eta_{T,k}^{\mathrm{term}}
 :=\calL_{T,k}+\frac12\log C_{T,d}^{\diag},
 \qquad k\in\mathcal T_T.
 \label{eq:terminal-eta}
\end{equation}
Finally, let \(\eta_{T,k}\) equal the bulk or terminal expression according
as \(k\in\mathcal B_T\) or \(k\in\mathcal T_T\).

\begin{proposition}[One-block continuation]
\label{prop:continuation}
The deterministic quantities \(\eta_{T,k}\ge0\) defined above are uniform
over all reachable histories and satisfy
\begin{equation*}
 \rho_{T,k}(\Xi_k)\le e^{\eta_{T,k}},\qquad
 \sum_{k<K_T}\eta_{T,k}=o(T).
\end{equation*}
\end{proposition}

\begin{proof}
Fix a reachable history $\Xi_k=(k,R_k,D_k),$
and write
\[
  q=q_k,
  \qquad
  q'=q_{k+1},
  \qquad
  A=\delta_T^d|D_k|,
  \qquad
  w=\Delta V_{T,k}(\Xi_k,Y),
  \qquad
  V=A+w.
\]
Split the next-guide numerator into its volume and escape branches, while
retaining the common full-guide denominator: $W_{T,k}
  =
  W_{T,k}^{\mathrm{vol}}
  +
  W_{T,k}^{\esc}.$

We first consider a bulk block \(q'\ge\theta_T\). Let
\(\mathcal D_{T,k}\) be the event that at least one endpoint in the fresh
candidate block receives the dagger token. On
\(\mathcal D_{T,k}^c\), Lemma~\ref{lem:continuation-estimates} gives
\[
  V\ge U_{T,k}^{\raw},
  \qquad
  \{V<v\}
  \subset
  \{U_{T,k}^{\raw}<v\}.
\]
Since the full current guide dominates its volume branch,
\begin{equation}
  W_{T,k}^{\mathrm{vol}}
  \le
  \exp\left\{
    -T\bigl[
      w+\JT_T(q',V)-\JT_T(q,A)
    \bigr]
  \right\}.
  \label{eq:volume-ratio}
\end{equation}
We bound the expectation of this quantity on
\(\mathcal D_{T,k}^c\) by four disjoint volume sectors.

If \(V<u_T\), then \(A\le V<u_T\), and for all sufficiently large \(T\)
both \(\JT_T(q,A)\) and \(\JT_T(q',V)\) lie on their linear branches.
Consequently, the right-hand side of \eqref{eq:volume-ratio} equals
\(e^{\calL_{T,k}}\). By
\eqref{eq:history-removal}, Proposition~\ref{prop:fresh-tail}, and
\eqref{eq:app-release-margin},
\[
  \begin{aligned}
  &\E\!\left[
    W_{T,k}^{\mathrm{vol}};
    V<u_T,\,
    \mathcal D_{T,k}^c
    \,\middle|\,
    \Xi_k
  \right]\le
  e^{\calL_{T,k}}
  \Pp\{U_{T,k}^{\raw}<u_T\mid\Xi_k\}\le
  \exp\left\{
    \calL_{T,k}
    -b_{T,k}\widetilde\lambda_T(u_T)
    +r_{T,k}
  \right\}
  \le
  e^{r_{T,k}}.
  \end{aligned}
\]

If \(u_T\le V\le v_T^{\up}\), then \(V\) takes at most
\(N_T+1\) grid-volume values. On a level \(V=v\),
Lemma~\ref{lem:scalar-continuation} and
\(T(q-q')=b_{T,k}\) imply $ W_{T,k}^{\mathrm{vol}}
  \le
  e^{b_{T,k}\widetilde\lambda_T(v)}.$
Moreover, on \(\mathcal D_{T,k}^c\), $\{V=v\}
  \subset
  \{U_{T,k}^{\raw}\le v\}.$
Therefore Proposition~\ref{prop:fresh-tail} gives 
\[
  \begin{aligned}
  \E\!\left[
    W_{T,k}^{\mathrm{vol}};
    V=v,\,
    \mathcal D_{T,k}^c
    \,\middle|\,
    \Xi_k
  \right]\le
  e^{b_{T,k}\widetilde\lambda_T(v)}
  \Pp\{U_{T,k}^{\raw}\le v\mid\Xi_k\}
  \le
  e^{r_{T,k}}.
  \end{aligned}
\]
Summing over the attainable levels yields
\[
  \E\!\left[
    W_{T,k}^{\mathrm{vol}};
    u_T\le V\le v_T^{\up},\,
    \mathcal D_{T,k}^c
    \,\middle|\,
    \Xi_k
  \right]
  \le
  (N_T+1)e^{r_{T,k}}.
\]

Next suppose that \(V>v_T^{\up}\) and
\(A\le A_T^{\cut}\). Since $v_T^{\up}
  =
  A_T^{\cut}+J+3,$
we have $w=V-A>J+3.$
Also, $H_{T,k+1}^{\mathrm{vol}}\le1$ and $ \JT_T(q,A)\le\Jc_Tq^p\le\Jc_T.$
Since \(\Jc_T\to J\), for all sufficiently large \(T\),
\[
  W_{T,k}^{\mathrm{vol}}
  \le
  \exp\{-T[w-\JT_T(q,A)]\}
  \le
  e^{-2T}.
\]

Finally, on the remaining sector with $V>v_T^{\up}$ and $A>A_T^{\cut},$
the scalar continuation inequality gives $W_{T,k}^{\mathrm{vol}}
  \le
  e^{b_{T,k}\widetilde\lambda_T(V)}.$
Because \(V\ge A>A_T^{\cut}\) and
\(\widetilde\lambda_T\) is decreasing, $W_{T,k}^{\mathrm{vol}}
  \le
  e^{b_{T,k}\widetilde\lambda_T(A_T^{\cut})}.$

Combining the four sectors,
\[
  \begin{aligned}
  \E\!\left[
    W_{T,k}^{\mathrm{vol}};
    \mathcal D_{T,k}^c
    \,\middle|\,
    \Xi_k
  \right]\le
  (N_T+2)e^{r_{T,k}}
  +
  e^{-2T}
  +
  e^{b_{T,k}\widetilde\lambda_T(A_T^{\cut})}.
  \end{aligned}
\]
Lemma~\ref{lem:continuation-estimates} also gives
\[
  \E\!\left[
    W_{T,k}^{\esc}
    \,\middle|\,
    \Xi_k
  \right]
  +
  \E\!\left[
    W_{T,k}^{\mathrm{vol}};
    \mathcal D_{T,k}
    \,\middle|\,
    \Xi_k
  \right]
  \le
  C_d
  e^{C_d\gamma_T^2b_{T,k}}.
\]
After increasing the fixed dimension-dependent constant \(C_d\) if
necessary, these bounds imply
\[
  \rho_{T,k}(\Xi_k)
  =
  \E[W_{T,k}\mid\Xi_k]
  \le
  \exp\{\eta_{T,k}^{\mathrm{bulk}}\},
\]
where \(\eta_{T,k}^{\mathrm{bulk}}\) is defined in
\eqref{eq:bulk-eta}.

Now consider a terminal block with \(q'<\theta_T\). The Cauchy--Schwarz
inequality and Lemma~\ref{lem:diagonal-envelope} give
\[
  \begin{aligned}
  \rho_{T,k}(\Xi_k)
  \le
  \mathfrak{d}_{T,k}(\Xi_k)^{1/2}\le
  e^{\calL_{T,k}}
  (C_{T,d}^{\diag})^{1/2}=
  \exp\{\eta_{T,k}^{\mathrm{term}}\},
  \end{aligned}
\]
where \(\eta_{T,k}^{\mathrm{term}}\) is defined in
\eqref{eq:terminal-eta}.

Thus \(\rho_{T,k}(\Xi_k)\le e^{\eta_{T,k}}\) uniformly over all reachable
histories. Finally, \eqref{eq:app-continuation-sum} gives $\sum_{k<K_T}\eta_{T,k}=o(T).$
For \(T\) below the deterministic asymptotic threshold,
the \(\eta_{T,k}\) may be enlarged without changing the conclusion.
\end{proof}

\subsection{Final assembly}
Finally, we prove the asymptotic optimality of the method.
\begin{theorem}[Natural-radius asymptotic optimality]
\label{thm:log-eff}
For every fixed \(d\ge2\) and \(a>0\), Algorithm~\ref{alg:package-w} is
asymptotically optimal in the sense of Definition~\ref{def:log-eff}.
\end{theorem}

\begin{proof}
From
\eqref{eq:pair-mass}, Lemma~\ref{lem:diagonal-envelope}, and
Proposition~\ref{prop:continuation},
\[
 \begin{aligned}
 m_{T,k}(\Xi_k)
 \le K_T^{-2}+e^{2\eta_{T,k}}\le e^{2\eta_{T,k}}(1+K_T^{-2}).
 \end{aligned}
\]
Since \(0\le\phi_T\le1\),
\[
 \begin{aligned}
 \log\E_{\Pi_T}\left[
   \phi_T(\Xi_{K_T})\prod_{k<K_T}m_{T,k}(\Xi_k)
 \right]\le2\sum_{k<K_T}\eta_{T,k}
 +K_T\log(1+K_T^{-2})=o(T).
 \end{aligned}
\]
Insert this estimate and Lemma~\ref{lem:initial-guide} into
\eqref{eq:exact-second-moment} to obtain $\log\E_{\Psi_T}Z_T^2\le-2JT+o(T).$
The matching lower bound follows from Theorem~\ref{thm:unbiasedness},
Proposition~\ref{prop:target-rate}, and Jensen's inequality:
\[
 \E_{\Psi_T}Z_T^2\ge(\E_{\Psi_T}Z_T)^2
 =\alpha_T^2=e^{-2JT+o(T)}.
\]
This finishes the proof.
\end{proof}

Since the relative variance of \(Z_T\) grows at most subexponentially
in \(T\), a subexponential number of i.i.d. replications suffices to
achieve any fixed relative variance bound or any fixed relative
accuracy at a prescribed confidence level.
Combining this sample requirement with the subexponential expected
work per replication in Proposition~\ref{prop:well-defined} gives
subexponential expected total work; Section~\ref{sec:numerics}
illustrates the tradeoff between sample size and replication cost.

\section{Finite-horizon lattice-surrogate experiments}
\label{sec:numerics}

We evaluate a finite-time lattice surrogate of the model in
Section~\ref{sec:setting} in six cases:
\[
(d,T)=(2,3),\ (2,3.5),\ (2,4),\ (5,0.5),\ (5,0.6),\ (5,0.7).
\]
We compare the statistical and computational performance of the lattice
analogue of Algorithm~\ref{alg:package-w} with direct Monte Carlo, and
illustrate the geometry of individual paths at $(d,T)=(2,4)$.
Appendix~\ref{sec:numericsdetail} gives the surrogate scaling and its relation
to the continuous construction.

\subsection{Surrogate target and evaluation metrics}
For each fixed pair $(d,T)$, let $S$ be a rate-one continuous-time Markov chain on $\mathbb Z^d$, started from the origin, with
uniform jump directions. Set
\[
 \eps=T^{-1/d},\qquad
 h=\eps/5,\qquad
 t_{\mathrm{raw}}=\frac{dT}{h^2},\qquad
 \beta=Th^d=5^{-d}.
\]
For $x\in\mathbb Z^d$, let
$B_{\mathbb Z^d}(x,5):=\{z\in\mathbb Z^d:\lVert z-x\rVert_2\leq5\}$,
and define the discrete sausage
$A:=\bigcup_{0\leq s\leq t_{\mathrm{raw}}}B_{\mathbb Z^d}(S_s,5)$.
The simulated target is
$\alpha^{\mathrm{lat}}_{d,T}:=\mathbb E[\exp\{-\beta|A|\}]$.
Thus Brownian motion is replaced by the lattice chain, while the continuous
Wiener sausage is replaced by the exact union of discrete Euclidean balls
centered at the visited sites. The choice $h=\eps/5$ represents the natural
sausage radius by a fixed lattice radius of five.

The lattice analogue of Algorithm~\ref{alg:package-w} scores independent
base-kernel candidate blocks using their exact new occupied volume and a
continuation guide, selects one candidate, and retains the corresponding
uniform-label/proposal-label likelihood correction. Within each stage,
conditionally independent candidate continuations can be generated and
scored in parallel with the GPU.

For independent replications $Y_1,\ldots,Y_n$ of an estimator $Y$, let
\[
 \overline Y=\frac{1}{n}\sum_{i=1}^nY_i,
 \qquad
 \widehat{\mathrm{CV}}=\frac{s_Y}{\overline Y},
 \qquad
 \operatorname{relCI}=\frac{1.96\,\widehat{\mathrm{CV}}}{\sqrt n},
\]
where $s_Y$ is the sample standard deviation and $\operatorname{relCI}$ is
the nominal relative $95\%$ confidence half-width. For each method and case,
we evaluate this half-width at prespecified sample-size checkpoints and stop
at the first checkpoint for which $\operatorname{relCI}\leq0.10$.
Checkpoint spacings differ across methods and cases, so this compares the
methods under a common stopping threshold; their realized relative
half-widths need not be equal.

\subsection{Algorithm~\ref{alg:package-w} versus direct Monte Carlo}

\begin{table}[htbp]
\caption{A100 GPU comparison of direct Monte Carlo and the
candidate-augmented lattice estimator. Each run stops at the first
method-specific checkpoint for which the estimated relative \(95\%\)
confidence half-width is at most \(10\%\).}
\label{tab:numerical-surrogate-results}
\vspace{4pt}
\centering
\small
\resizebox{\textwidth}{!}{%
\begin{tabular}{c l r c r c r r}
\hline
$(d,T)$
& method
& $n$
& estimate
& CV
& relative half-width
& GPU time (s)
& ms/rep \\
\hline
$(2,3)$
& Direct MC
& 70,000
& $1.44\times10^{-9}$
& 13.3503
& 9.89\%
& 4.1
& 0.0586 \\
$(2,3)$
& Augmented IS
& 600
& $1.32\times10^{-9}$
& 1.1376
& 9.10\%
& 1.4
& 2.373 \\
\hline
$(2,3.5)$
& Direct MC
& 280,000
& $7.09\times10^{-11}$
& 25.9446
& 9.61\%
& 18.0
& 0.0642 \\
$(2,3.5)$
& Augmented IS
& 1,000
& $6.71\times10^{-11}$
& 1.5996
& 9.91\%
& 3.9
& 3.862 \\
\hline
$(2,4)$
& Direct MC
& 950,000
& $3.12\times10^{-12}$
& 49.2644
& 9.91\%
& 65.8
& 0.0692 \\
$(2,4)$
& Augmented IS
& 2,200
& $3.10\times10^{-12}$
& 2.3688
& 9.90\%
& 8.6
& 3.921 \\
\hline
$(5,0.5)$
& Direct MC
& 30,000
& $1.52\times10^{-9}$
& 8.7839
& 9.94\%
& 20.9
& 0.696 \\
$(5,0.5)$
& Augmented IS
& 800
& $1.48\times10^{-9}$
& 1.2369
& 8.57\%
& 13.7
& 17.160 \\
\hline
$(5,0.6)$
& Direct MC
& 140,000
& $7.74\times10^{-11}$
& 17.4838
& 9.16\%
& 124.9
& 0.892 \\
$(5,0.6)$
& Augmented IS
& 800
& $7.92\times10^{-11}$
& 1.4422
& 9.99\%
& 17.7
& 22.186 \\
\hline
$(5,0.7)$
& Direct MC
& 360,000
& $4.12\times10^{-12}$
& 30.2444
& 9.88\%
& 397.2
& 1.103 \\
$(5,0.7)$
& Augmented IS
& 2,600
& $4.63\times10^{-12}$
& 2.5693
& 9.88\%
& 52.0
& 20.001 \\
\hline
\end{tabular}%
}
\end{table}

Table~\ref{tab:numerical-surrogate-results} compares direct Monte Carlo and
the candidate-augmented lattice estimator on the same NVIDIA
A100-SXM4-40GB GPU. In all six cases, the two estimates agree within their
combined Monte Carlo uncertainty, and all reported nominal relative $95\%$
half-widths are at most $10\%$. The augmented estimator reduces the observed
coefficient of variation by factors ranging from $7.10$ to $20.80$, and
reduces the required sample size by factors ranging from $37.5$ to $431.8$.
Despite its higher cost per replication, it reduces the observed computational
time by factors of approximately $1.5$--$7.6$.

The improvement becomes more pronounced as the two-dimensional target
becomes rarer. As $T$ increases from $3$ to $4$, the coefficient-of-variation
reduction increases from $11.74$ to $20.80$, while the observed GPU time
reduction increases from approximately $2.9$ to $7.6$. At $(d,T)=(2,4)$,
direct Monte Carlo uses $950{,}000$ replications and $65.8$ seconds, whereas
the augmented estimator uses $2{,}200$ replications and $8.6$ seconds.

The three five-dimensional cases show a similar finite-horizon pattern.
The observed GPU time reduction is approximately $1.5$, $7.0$, and
$7.6$ at $T=0.5$, $0.6$, and $0.7$, respectively. In the rarest case,
$(d,T)=(5,0.7)$, the two methods estimate a target of order $10^{-12}$;
direct Monte Carlo uses $360{,}000$ replications and $397.2$ seconds, whereas
the augmented estimator uses $2{,}600$ replications and $52.0$ seconds.

The augmented estimator is more expensive per replication in every reported
case, as expected from its candidate-generation and selection steps. At
$(d,T)=(5,0.5)$, for example, direct Monte Carlo costs approximately $0.696$
milliseconds per replication, compared with $17.160$ milliseconds for the
augmented estimator. Its lower variability compensates for this additional
cost and yields a smaller total computational time.

\subsection{Path geometry and repeated visits}

\begin{figure}[H]
\centering
\includegraphics[width=\textwidth]{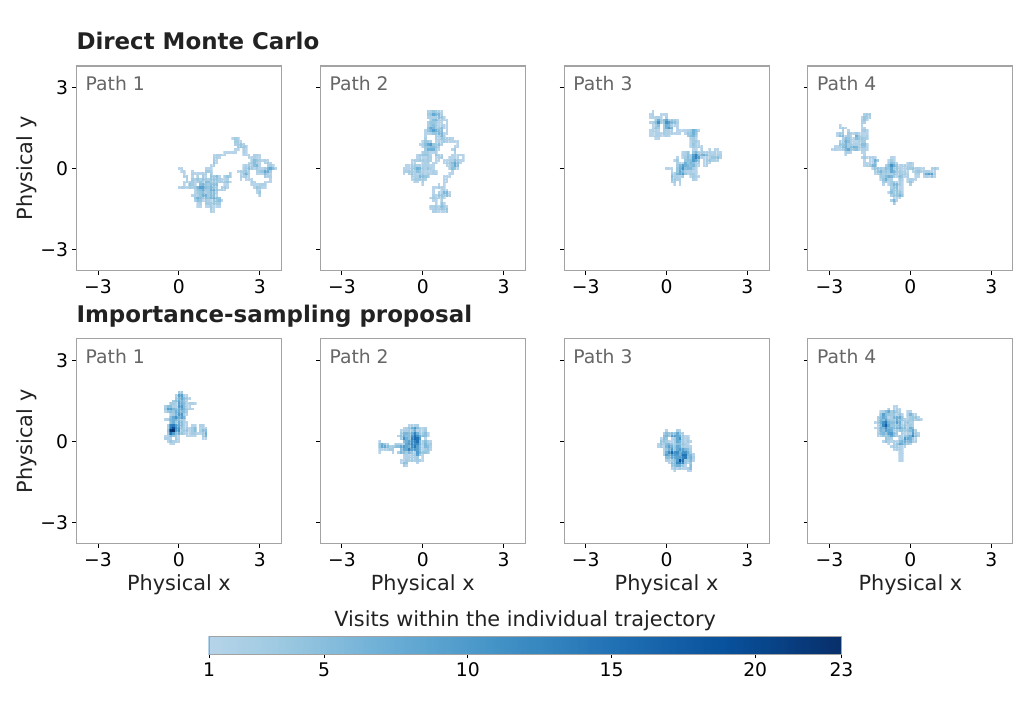}
\caption{Four direct-Monte-Carlo paths (top) and four importance-sampling proposal paths (bottom) at $(d,T)=(2,4)$. Paths are
selected uniformly without replacement.
All panels use the same physical axes and linear color scale. Each panel
shows an individual lattice trajectory; colors indicate the number of
visits to each site, including the initial origin and subsequent returns.}
\label{fig:path-geometry}
\end{figure}

Figure~\ref{fig:path-geometry} shows four uniformly selected trajectories
from each method at $(d,T)=(2,4)$, using the same configuration as
Table~\ref{tab:numerical-surrogate-results}. In these examples, the
importance-sampling proposal explores more compact, often approximately
ball-like regions and revisits sites more frequently. This is consistent
with the volume--confinement intuition behind the guide: the
Donsker--Varadhan variational balance trades explored volume against the
spectral cost of confinement, and the Faber--Krahn inequality identifies a
ball as minimizing that cost at fixed volume. The figure provides a
qualitative finite-horizon illustration of this intuition.

\section{Discussion and future directions}
\label{sec:discussion}
We develop an exact sequential importance sampler for the survival
probability of a Brownian particle among fixed-radius Poissonian
obstacles, using its equivalent formulation under Brownian scaling
as a natural-radius Wiener-sausage Laplace transform. For every fixed $d\geq2$
and $a>0$, Algorithm~\ref{alg:package-w} is unbiased at every finite horizon
by Theorem~\ref{thm:unbiasedness} and asymptotically optimal as $T\to\infty$
by Theorem~\ref{thm:log-eff}. Proposition~\ref{prop:well-defined} also
establishes finite expected total work at each fixed horizon and
subexponential expected total work as \(T\to\infty\). 
The construction replaces
repeated path-dependent spectral calculations on irregular explored regions
by conditionally independent candidate clouds, score-based selection, and an
exact proper-weight identity.  The inward grid certificate and residual
Poisson correction recover the continuous target,
while the candidate-cloud structure permits the exact
diagonal--off-diagonal second-moment representation in
Theorem~\ref{thm:exact-second} and is naturally suited to parallel
implementation.
The present proof is specific to standard Brownian motion and the
natural-radius scaling $\eps_T=aT^{-1/d}$.  A natural theoretical direction
is to extend the construction to more general diffusions, including models
with drift, spatially varying covariance, or heterogeneous environments, as
well as to other swept-region functionals and sausage geometries.  The
candidate-augmentation and proper-weight mechanisms should remain applicable,
but such extensions would require model-specific analogues of the confinement
guide, the fresh-block lower-volume estimate, and the exact residual
correction.

As discussed in the Introduction, swept-volume functionals arise in molecular
communication with absorbing receivers, diffusion-limited chemical and
biological reactions, finite-sized particle transport in porous media, and
the volume explored by mobile excitations in materials.  In these settings,
quantities of interest may involve rare survival, no-detection, or encounter
probabilities determined by the region swept out by a diffusing trajectory.
The methodology developed here may therefore provide a useful basis for
efficient simulation when direct Monte Carlo is dominated by atypically
compact paths and consequently has large relative variance.  Applying the
method to any particular physical model would require incorporating its
specific dynamics, interaction geometry, and observation mechanism, together
with a corresponding second-moment analysis; these application-specific
extensions are left for future work.

\appendix

\section{Deterministic certification and Brownian-bridge completion}
\label{app:exact-real}

\subsection{Deterministic and kernel constructions}
\label{app:deterministic-kernels}
\label{app:constants}

We construct the constants in \eqref{eq:certified-constants} by rational
enclosures, with all operations charged as in Appendix~\ref{app:work-model}.
The constants in \eqref{eq:certified-constants} are determined by
the principal Dirichlet eigenvalue \(\lambda_d\). This eigenvalue satisfies
\(\lambda_d=j^2\omega_d^{2/d}\), where \(\omega_d\) is the volume
of the unit-radius ball and \(j\) is the smallest positive zero
of the Bessel function \(J_\nu\) of the first kind,
with \(\nu=d/2-1\). To approximate this zero with a rigorous error bound, we use the normalized series
\[
 F_d(z)=\sum_{m=0}^{\infty}
 \frac{(-z^2/4)^m}{m!(d/2)_m},
 \qquad
 (d/2)_m=\prod_{i=0}^{m-1}(d/2+i),
\]
with the empty product equal to one.
For \(z>0\),
\(J_\nu(z)=(z/2)^\nu F_d(z)/\Gamma(\nu+1)\), so the positive
prefactor implies that \(F_d\) and \(J_\nu\) have the same positive
zeros and signs. The series has rational terms at rational arguments,
which allows certified evaluation using rational arithmetic.
Let \(j\) denote their first positive zero.
This zero is simple and satisfies \(F'_d(j)<0\). We construct an initial enclosure of \(j\) by a deterministic search.
Each candidate specifies rational numbers \(0<a<b\) and \(m_d>0\),
a finite rational subdivision of \([0,a]\), and series truncation
orders. We examine the candidates in increasing order of their total
integer encoding size, using lexicographic order to break ties.
We accept the first candidate whose interval bounds certify that
\(F_d>0\) on \([0,a]\), \(F_d(a)>0>F_d(b)\), and
\(F'_d\le-m_d\) on \([a,b]\).

Such a finite certificate exists because \(F_d\) has a strictly
positive minimum on each compact interval before \(j\), while
\(F'_d(j)<0\) and continuity ensure that the derivative remains
bounded away from zero on a sufficiently small interval around \(j\).
We verify these inequalities using polynomial interval evaluation
and uniform bounds on the series remainders
\cite{MooreKearfottCloud2009Interval}.

To make these remainder bounds explicit, we fix a positive rational
number \(B\ge b\) and bound the absolute terms of both series
uniformly on \(|z|\le B\).
The successive ratios of these bounds are at most
\(B^2/[4(m+1)^2]\) for \(F_d\) and \(B^2/(4m^2)\) for \(F'_d\),
with \(m\ge1\) in the latter case.
Once these ratios are at most \(1/2\), twice the next term bounds the
tail. The initial search therefore has finite cost depending only on \(d\).

At each refinement, we evaluate \(F_d\) at the midpoint \(x\)
of the current root enclosure, of width \(w\), with a certified
value interval of width at most \(m_dw/8\).
A certified sign determines which half to retain.
Otherwise, the derivative bound gives $|x-j|\le \frac{|F_d(x)|}{m_d}\le \frac{w}{8},$
so we retain \([x-w/8,x+w/8]\).
Each step therefore preserves the root enclosure and reduces
its width to \(w/2\) or \(w/4\).
Starting from a fixed positive rational enclosure of \(\pi\),
we refine it to width at most \(w\) using the alternating series
in Machin's identity
\(\pi=16\arctan(1/5)-4\arctan(1/239)\).
We then propagate the enclosures of \(j\) and \(\pi\) through
\[
 \omega_0=1,\quad \omega_1=2,\quad
 \omega_m=\frac{2\pi}{m}\omega_{m-2}\ (m\ge2),\qquad
 \lambda_d=j^2\omega_d^{2/d},\qquad
 v_\star^{d+2}=\frac{j^{2d}\omega_d^2}{d^d}.
\]
The last identity gives positive rational bounds
\(q_-\le v_\star^{d+2}\le q_+\).
Using bisection with exact integer-power comparisons, we enclose
\(q_\pm^{1/(d+2)}\) to width at most \(w\) and obtain
\(v_\star\in[u,v]\) with \(v-u=O_d(w)\).
We use fixed positive rational brackets, obtained initially by
doubling and halving, and denote their common upper bound by \(V\).
We stop when
\[
 \frac{v-u}{2}\left[1+\frac{d+2}{2}
                    +\frac{d+2}{2}\max\{1,V\}\right]
 \le(1\vee T)^{-4}.
\]
Then we set \(\widetilde v_{\star,T}=(u+v)/2\) and define
\(\Jc_T\) and \(\chic_{d,T}\) by the coupling identities in
\eqref{eq:certified-constants}, evaluating the fractional power
as \(\exp\{(1+2/d)\log\widetilde v_{\star,T}\}\).
The midpoint error is at most \((v-u)/2\).
Together with the linear coupling for \(\Jc_T\) and the derivative
bound
\[
 \frac{d+2}{2}y^{2/d}
 \le \frac{d+2}{2}\max\{1,V\},
 \qquad 0<y\le V,
\]
this proves that the stopping rule controls the sum of the three
errors in \eqref{eq:certified-constants}.

After \(s\) refinements, the root-bracket width lies between
\(4^{-s}w_0\) and \(2^{-s}w_0\), where \(w_0\) is its initial
width. The geometric remainder bounds and recursive series
evaluation, together with the endpoint-root bisections, therefore
give a cost of \(O_d(s+1)\) operations per stage.
Since \(v-u=O_d(w)\), only \(O_d(1+\log(1\vee T))\) stages
are needed. Including the finite initialization cost, we obtain
\begin{equation}
 C_{\rm cert}(T)=O_d\!\left((1+\log(1\vee T))^2\right).
 \label{eq:app-certification-cost}
\end{equation}

The half-open cells and dagger convention make the token map Borel;
the dagger retains the real endpoint and adds no paid cell.
Since \(1\le K_T\le N_T\), \(n_{T,k}\ge1\),
\(1\le M_{T,k}<\infty\), and every score in \eqref{eq:scores} is positive,
the Gaussian and categorical kernels define \(\Psi_T^{\sk}\) by
Ionescu--Tulcea. Moreover,
\begin{equation}
 d_T^{\cem}(0)=L_T^{\rm win}\delta_T\ge T^2-\delta_T.
 \label{eq:app-cem-origin}
\end{equation}
The inclusion \(\overline Q_{T,m}\subset\overline B_2(x,\eps_T)\)
for \(x\in Q_{T,m}\), including cell faces, follows from
\(\delta_T\sqrt d<\eps_T\) and proves the paid-set inclusion.

\subsection{Initialization of persistent Brownian bridges}
\label{app:bridge-init}

For each coordinate of each fine gap, we use the
L\'evy--Ciesielski representation of the centered bridge on \([0,h]\),
with height-one dyadic tents
\cite{BierkensGrazziVanDerMeulenSchauer2021Bridges,
BrownGriebelKuoSloan2024LevyCiesielski}:
\[
 \beta(t)=\sum_{n\ge0}\sum_{0\le r<2^n}
 \sqrt h\,2^{-(n+2)/2}Z_{n,r}\psi_{n,r}(t),\qquad
 Z_{n,r}\stackrel{\rm iid}{\sim}N(0,1),\quad h=h_T.
\]
Here \(\psi_{n,r}\) is supported on
\([rh2^{-n},(r+1)h2^{-n}]\). With the standard normal random variable $Z\sim N(0,1)$, we set
\[
 b_n=\sqrt{8(n+1)},\quad c_n=2\,32^{-(n+1)},\quad
 \lambda_n=-\log(1-c_n),\quad
 \vartheta_n=\Pp\{|Z|>b_n\}.
\]
We call a coefficient exceptional if \(|Z_{n,r}|>b_n\).

Each coordinate of each fine gap has its own coefficient array
\((Z_{n,r})_{n\ge0,\,0\le r<2^n}\).
We initialize these arrays independently. We draw
\(H\sim\operatorname{Poisson}(2/15)\) labels, sampling each level
\(n\) with probability \((15/16)16^{-n}\) and each index \(r\)
uniformly from \(\{0,\ldots,2^n-1\}\).
We retain each label with probability \(\lambda_n/(2c_n)\)
and collapse duplicate keys.
For each retained key, we draw one independent proposal
\(G\sim N(0,1)\), \(U\sim\operatorname{Unif}(0,1)\), accepting when both $G\ge0$ and $U\le(32e^{-4})^{n+1}e^{-b_nG}.$
On acceptance, we mark the key and store \(Z_{n,r}=S(b_n+G)\),
with an independent fair sign \(S\).
For each unmarked key, we draw fresh standard normals when first
requested until \(|Z|\le b_n\).
We discard failed initializer proposals without using them as coefficients.
For each array, we retain all marks and exposed coefficients so that
subsequent residual-membership queries use the same bridge realization.

The label count at each key has mean \(2c_n\), and
\(c_n\le\lambda_n<2c_n\). Thinning and duplicate collapse therefore
give independent potential-key indicators of probability \(c_n\).
A direct Gaussian density calculation gives acceptance probability
\(\vartheta_n/c_n\) for each retained key.
The accepted value has the outside-truncated normal law, and the
threshold is at most one for \(G\ge0\), since \(32e^{-4}<1\).
Thus the independent marks and conditional inside/outside laws
reconstruct the iid normal array.
Each inside rejection loop has mean at most \(16/15\).
Across the \(D=dN_T\) arrays, we let \(L_{\max}\) denote the largest
marked level, or \(-1\) if none is marked. Since \(\vartheta_n\le c_n\), then
\begin{equation}
 \Pp\{L_{\max}\ge n\}\le\min\{1,D16^{-n}/15\},\qquad n\ge0.
 \label{eq:app-exception-tail}
\end{equation}
We add the selected endpoint interpolants to obtain conditional
Brownian bridges, independently across gaps and coordinates.

We write \(P_q\) for the scalar truncation through level \(q\).
If all higher levels are unmarked, then
\begin{equation}
 \|\beta-P_q\|_\infty
 \le\sqrt h\sum_{n>q}2^{-(n+2)/2}b_n
 \le\rho_q(h):=16\sqrt h\sqrt{q+2}\,2^{-q/2},\qquad q\ge0.
 \label{eq:app-bridge-tail}
\end{equation}
The last inequality follows from
\(\sum_{j\ge1}\sqrt j\,r^j\le r/(1-r)^2\), with \(r=2^{-1/2}\).
We define
\[
 e_q=\sqrt d\,\rho_q(h_T),\quad
 q_*=\min\{q\ge0:e_q\le\eps_T/12\},\quad
 q_0=\max\{q_*,L_{\max}\}.
\]
We expose all arrays through \(q_0\), storing binary trees of
midpoint values and child pointers.
We let \(A\) be the maximum \(\ell_\infty\)-norm of the polygon
vertices after adding the endpoint interpolants, and set
\[
 \mathcal R_T=1+\left\lceil A+\rho_{q_0}(h_T)\right\rceil,
 \qquad\mathcal C_T^{\res}
 =[-\mathcal R_T-\eps_T,\mathcal R_T+\eps_T)^d.
\]
This box contains the entire sausage.
Later tents vanish at existing vertices, so
\(A\le S_T:=\sup_{t\le T}\|B_t\|_\infty\). Hence $\mathcal R_T\le S_T+\rho_0(h_T)+2.$

\subsection{A finite residual Poisson query}
\label{app:residual-query}
After constructing the box, we use fresh randomness to draw a Poisson
count of mean \(T|\mathcal C_T^{\res}|\) and that many independent
uniform points in it.
Conditional on the path and initialization data, these points form
an intensity-\(T\) Poisson process on \(\mathcal C_T^{\res}\).
We skip points in \(\calI_T\) and process the fine gaps in a fixed
order for each remaining point \(z\).
In each gap, we initialize a first-in-first-out (FIFO) queue with
the \(2^{q_0+1}\) initial polygon intervals in time order.
For an interval at truncation level \(q\), let \(s\) be the distance
from \(z\) to its polygon segment.
We return a hit if \(s<\eps_T-e_q\), prune the interval if
\(s>\eps_T+e_q\), and otherwise expose its midpoint in every
coordinate and append its two children, left first.
New coefficients have level \(q+1>q_0\) and use the inside-truncated
sampler. For each point, we restart the same FIFO traversal from
the initial roots and reuse each array's stored coefficient values.
The local error is at most \(e_q\) by \eqref{eq:app-bridge-tail}.

A hit in any gap proves membership; empty queues in every gap prove
nonmembership.
We return \(Y_T^{\res}=0\) at the first hit outside \(\calI_T\),
and \(Y_T^{\res}=1\) if all point queries finish without a hit.
Uniform convergence and FIFO traversal ensure that each gap query
terminates unless \(z\) is exactly at distance \(\eps_T\) from
that gap's path range.
For every nonempty compact \(K\) and \(r>0\), the level set
\(\{z:\dist(z,K)=r\}\) is Lebesgue-null: the distance function has
gradient of norm one almost everywhere outside \(K\), and gradient
zero almost everywhere on each level set.
The finite union of these gapwise tangency sets is therefore
Lebesgue-null and is almost surely avoided by the Poisson list.
Half-open cell comparisons and deterministic equality conventions
give measurable tests.
We assign a fixed output on the remaining null event to extend
the kernel everywhere.
Because the box contains the sausage and was chosen without Poisson
randomness, then
\[
 \Pp\{Y_T^{\res}=1\mid B,\text{augmented grid data}\}
 =\exp\{-T|W_T^{\eps_T}\setminus\calI_T|\}.
\]

\subsection{Operation model and integer operations}
\label{app:work-model}

We use an exact-real scalar-operation model of the type described in
\cite[Section~2.3 and Appendix~A]{BlanchetJinLiu2026Exact}.
The inputs \(d,a,T\), \(\pi\), and rational constants are exact.
Scalar arithmetic, comparison, square root, exponential, logarithm,
independent standard-normal or uniform generation, and elementary memory,
indexing, and output operations have fixed positive finite cost.
Vector operations are charged coordinatewise; integer arithmetic follows
the same convention, without charging bit length. Floors and ceilings
are implemented by counted comparison-and-increment loops.
Let \(\mathsf{Work}_T\) include all these operations in one run of
Algorithm~\ref{alg:package-w}, including certification, rejected draws,
initialization, and repeated tests of cached nodes.

We generate a Poisson count of mean \(\lambda\) using exponential
interarrivals \(-\log U\), with expected cost \(O(1+\lambda)\).
Uniform trials with success probability \(15/16\), followed by
\(n\) fair bits, generate an initializer level \(n\) and its index
in \(O(n+1)\) operations.
Window checks and linear scans give tokenization and paid-cell
membership costs \(O_d(1+L_T^{\rm win})\) and \(O_d(1+N_T)\),
respectively.
Counted loops implement candidate and enclosure ceilings at costs
\(O(1+M_{T,k})\) and \(O(1+\mathcal R_T)\); other mesh and block
loops have polynomial cost in \(N_T+L_T^{\rm win}\).
We store each node's depth, scale, endpoints, and child pointers.
Midpoint and error-bound updates, segment-distance calculations
(with a separate zero-length case), and queue operations each
cost \(O_d(1)\), with coefficient-sampling trials counted separately.

\subsection{Finite and subexponential expected total work}
\label{app:work-bound}

\begin{theorem}[Expected total work]
\label{thm:expected-work}
Fix \(d\ge2\) and \(a>0\).
Under the implementation and cost model above,
\(\E_{\Psi_T}\mathsf{Work}_T<\infty\) for every \(T>0\).
Set \(M_*:=\max_{k<K_T}M_{T,k}\).
There exist constants \(C,c,T_0>0\), depending only on \(d,a\), such that
\[
 \E_{\Psi_T}\mathsf{Work}_T\le C(1+T)^cM_*=e^{o(T)},
 \qquad T\ge T_0.
\]
\end{theorem}

\begin{proof}
For the convenience of the notations, we write \(N=N_T\), \(h=h_T\), \(\epsilon=\eps_T\), \(D=dN\), and
\(p=d/(d+2)\). Unless stated otherwise, expectations are under
\(\Psi_T\), and constants may depend on fixed \(d,a\).
For \(T\ge1\), there exist constants \(C,c>0\) such that
\[
 \frac{\sqrt h}{\epsilon}\le
 \frac1{8d\sqrt T\,\ell_T^{3/2}},\qquad
 N+L_T^{\rm win}\le C(1+T)^c.
\]
By \eqref{eq:certified-constants}, \(\Jc_T=O(1)\), so
\(\gamma_T=O(T^{-1})\). Since \(b_{T,\max}\le T\), we also have
\(C_{T,d}^{\diag}=O(1)\).
Using \(x^p-y^p\le(x-y)^p\) for \(x\ge y\ge0\) and
\eqref{eq:schedule-orders}, we obtain
\[
 \begin{gathered}
 \max_k\calL_{T,k}\le\Jc_T T^{1-p}b_{T,\max}^p=o(T),\\
 M_*\le1+C_{T,d}^{\diag}K_T^2e^{2\max_k\calL_{T,k}}=e^{o(T)},
 \qquad\sum_kM_{T,k}n_{T,k}\le NM_*.
 \end{gathered}
\]

\paragraph*{Local refinement}
We first establish a local bound on the expected refinement cost.
For a continuous gap \(w:[0,h]\to\R^d\), set
\(\Delta_q=h2^{-q}\), \(t_{q,j}=j\Delta_q\),
\(S_{q,j}(w)=[w(t_{q,j}),w(t_{q,j}+\Delta_q)]\), and
\(r_q=\epsilon-3e_{q-1}\). We define the functional
\begin{equation}
 \mathfrak A_T(w)=2\sum_{q\ge q_*+1}\sum_{0\le j<2^q}
 \int_{\R^d}\ind\{\dist(z,w[0,h])>r_q\}
 \ind\{\dist(z,S_{q,j}(w))\le\epsilon+e_{q-1}\}\,\dd z.
 \label{eq:app-work-envelope}
\end{equation}
Let \(B^{(i)}\) be the completed \(i\)-th gap, parametrized on
\([0,h]\), and let \(J_i(B,z)\) count the number of noninitial tests in a
standalone query that stops at its first hit.
A child tested at level \(q\) has an unresolved parent at level
\(q-1\), whose interval has length \(\Delta_q\).
FIFO order ensures that all surviving parents at that level have
already been tested without a hit. Their path pieces therefore have
distance at least \(\epsilon-2e_{q-1}>r_q\) from \(z\), while
previously pruned pieces have distance greater than \(\epsilon\).
Each parent has at most two children, so
\begin{equation}
 \int_{\R^d}J_i(B,z)\,\dd z\le\mathfrak A_T(B^{(i)}).
 \label{eq:app-work-pathwise}
\end{equation}

We next bound \(\E_{\rm BM}\mathfrak A_T(B)\), where \(B\) is a
standard Brownian motion on \([0,h]\).
Let \(S_r(s,\rho)\) be the probability of avoiding a radius-\(r\)
ball for time \(s\), starting at distance \(\rho\) from its center.
For \(r\ge\epsilon/2\), \(x>0\), and \(0<s\le h\),
\[
 S_r(s,r+x)\le\min\{1,C_Tx/\sqrt s\},\qquad
 C_T=\sqrt{2/\pi}+2(d-1)\sqrt h/\epsilon.
\]
Indeed, the radial drift before hitting the ball is at most
\((d-1)/(2r)\). Comparison with a one-dimensional Brownian motion having this
constant drift and the reflection formula for its running minimum
give \(\sqrt{2/\pi}\,x/\sqrt s+(d-1)x/r\).
The mesh bound makes \(C_T\) uniformly bounded for \(T\ge1\).
For \(\epsilon/2\le r\le\epsilon\), \(0<u\le\epsilon\), and
\(0<t<h\), spatial translation by \(B_t\) and independence of the
reversed past and future increments give
\begin{align*}
 &\E_{\rm BM}\int_{\R^d}
 \ind\{\dist(z,B[0,h])>r,\ r<|z-B_t|\le r+u\}\,\dd z\\
 &\quad=\int_{r<|y|\le r+u}S_r(t,|y|)S_r(h-t,|y|)\,\dd y
 \le\frac{C_d\epsilon^{d-1}u^3}{\sqrt{t(h-t)}}.
\end{align*}
At \(t=0,h\), the single survival factor instead gives
\(C_d\epsilon^{d-1}u^2/\sqrt h\).

Set \(s_q=4\sqrt d\sqrt{(q+1)\Delta_q}\) and
\(u_q=4e_{q-1}+s_q\). For \(q\ge q_*+1\),
\[
 r_q\ge3\epsilon/4,\qquad
 s_q=e_{q-1}/(4\sqrt2),\qquad u_q<\epsilon.
\]
When a parent increment has norm at most \(s_q\) (small increments), the integrand in
\eqref{eq:app-work-envelope} restricts \(z\) to
\(r_q<|z-B_{t_{q,j}}|\le r_q+u_q\).
The shell bounds and
\(\sum_{j=1}^{2^q-1}[j\Delta_q(h-j\Delta_q)]^{-1/2}\le4/\Delta_q\)
therefore bound the expected level contribution by $C_d\epsilon^{d-1}\sqrt h
 \bigl((q+1)^{3/2}2^{-q/2}+(q+1)2^{-q}\bigr).$
For larger increments, we obtain an upper bound by replacing
the avoidance indicator with one and retaining the second indicator.
The remaining capsule has volume at most
\(C_d(\epsilon^d+\epsilon^{d-1}|\Delta B|)\).
Standard Gaussian integral computations and Cauchy--Schwarz inequality bound this level contribution by
\(C_d2^q(\epsilon^d+\epsilon^{d-1}\sqrt h)e^{-4(q+1)}\).
Both series converge, thus
\begin{equation}
 \E_{\rm BM}\mathfrak A_T(B)
 \le C_d(\epsilon^d+\epsilon^{d-1}\sqrt h)
 \le C_d\epsilon^d,\qquad T\ge1.
 \label{eq:app-work-envelope-mean}
\end{equation}

\paragraph*{Initialization and enclosure}
We next control the initialization cost and the size of the enclosing box.
The centered bridge arrays and their initializer randomness are
independent of the selected skeleton.
For \(T\ge1\), the mesh bound \(q_*\le12\) and
\eqref{eq:app-exception-tail} imply
\(\E4^{q_0}\le C_d(1+D)\).
The label count \(H\sim\operatorname{Poisson}(2D/15)\) and geometric
level tails give polynomial expected cost for label generation and
duplicate removal.
Including exception-list scans and rejection attempts, the expected
cost of initial exposure is at most
\[
 C_dD\E[(1+H)2^{q_0}]
 \le C_dD(1+D)^{3/2},
\]
by Cauchy--Schwarz inequality.

The candidate increments use
\(N_T^{\mathrm G}=d\sum_kM_{T,k}n_{T,k}\le dNM_*\)
independent scalar normals. Therefore,
\[
 \E S_T^{2d}
 \le C_d\left[
 (N\sqrt h)^{2d}\bigl(1+\log(2N_T^{\mathrm G})\bigr)^d
 +Nh^d
 \right].
\]
Since \(\log(2N_T^{\mathrm G})=O(T)\) and
\(\mathcal R_T\le S_T+\rho_0(h)+2\),
Cauchy--Schwarz inequality gives
\[
 \E[2^{q_0}(\mathcal R_T+\epsilon)^d]
 \le
 \{\E4^{q_0}\,\E(\mathcal R_T+\epsilon)^{2d}\}^{1/2}
 \le C(1+T)^c.
\]

\paragraph*{Selection and total work}
We now transfer the local estimates to the selected path and combine
the work bounds.
Let \(\mathcal H_k\) denote the history before cloud \(k\).
For every nonnegative and measurable function \(f\), we have
\[
 \E[f(Y^{I_k})\mid\mathcal H_k]
 \le M_{T,k}\int f(y)\,p_{T,k}(\mathcal H_k,\dd y).
\]
This comparison extends to completed blocks by averaging the
conditional bridges. Translation invariance makes the reference
expectation of \(\mathfrak A_T\) independent of the incoming position.

Conditional on the grid, initialization data, and full bridge
coefficient arrays, we apply Campbell's formula
\cite[Proposition~2.7]{LastPenrose2017Poisson}
to standalone queries. Ignoring paid-cell skipping and early stopping, the block
comparison, \eqref{eq:app-work-pathwise}, and
\eqref{eq:app-work-envelope-mean} give
\[
 \E[\text{noninitial tests}]
 \le T\sum_{i=1}^N\E\mathfrak A_T(B^{(i)})
 \le C_dT\epsilon^d\sum_kM_{T,k}n_{T,k}
 \le C_da^dNM_*.
\]
Each gap incurs only its own block's selection factor.
Poisson generation, paid-cell scans, initial-root tests, and queue
initialization cost at most
\(C_d(1+NT\E[2^{q_0}|\mathcal C_T^{\res}|])\)
in expectation, which is polynomial by the preceding mixed moment bound.

Each unresolved parent requests at most \(d\) new coefficients
and creates at most two children, so node, cache, and queue
operations, including cleanup, are charged to initial nodes and
parent tests.
Fresh rejection streams start only upon request and have conditional
expected trial count at most \(16/15\); Tonelli's theorem therefore
bounds their total expected cost by a constant multiple of the
expected number of requests.

Candidate generation and processing, including cloud scans and
counted ceilings, cost at most \(NM_*\) times a polynomial in
\(N+L_T^{\rm win}\).
The remaining integer loops, enclosure ceiling, and certification
have polynomial expected cost in \(T\), by the preceding bounds and
\eqref{eq:app-certification-cost}.
Combining these estimates proves the stated bound.
For every fixed \(T>0\), all deterministic parameters are finite
and the same level sums and moment bounds remain finite, thus the proof is finished.
\end{proof}

\section{Technical details}
\label{app:lower-volume}

Throughout this section, \(C_d\in(0,\infty)\) denotes a finite constant
depending only on \(d\), whose value may change from line to line.
This appendix contains the technical details used in the proofs in the main
text. We first check the asymptotic schedule used in
Section~\ref{sec:efficiency}, which follows directly from the definitions.
Eventually,
\[
 \delta_T=\frac{a}{8d}T^{-1/d},\quad
 N_T=\Theta_{d,a}(T^{2+2/d}(\log T)^3),\quad
 h_T=\Theta_{d,a}(T^{-1-2/d}(\log T)^{-3}).
\]
For \(d=2\),
\[
 \mathfrak H_T=\Theta(\log T),\quad
 b_T=\Theta(\sqrt{T\log T}),\quad
 K_T=\Theta(\sqrt{T/\log T}),
\]
whereas, for fixed \(d\ge3\),
\[
 \mathfrak H_T=\Theta_{d,a}(T^{(d-2)/d}),\quad
 b_T=\Theta_{d,a}(T^{(d-1)/d}),\quad
 K_T=\Theta_{d,a}(T^{1/d}).
\]
Since \(|b_{T,k}-T/K_T|\le h_T\), the balanced block lengths are comparable.
These formulas prove \eqref{eq:schedule-orders} and give
\[
 b_T\delta_T^{d-2}\to\infty\quad(d\ge3),\qquad
 \frac{b_T}{\log(b_{T,\max}\delta_T^{-2})}\to\infty\quad(d=2).
\]
We now state the lemmas deferred from the main text.
\begin{lemma}[Grid filling and fixed-radius translation]
\label{lem:raster-source-translation}
Fix a stage \(k<K_T\), an incoming history \(\Xi_k\), and
\(\kappa,\zeta>0\). Conditional on the fresh Gaussian candidate endpoints
\[
  Y=(Y_1,\ldots,Y_{n_{T,k}}),
  \qquad
  Y_0=x(\Xi_k),
\]
complete each mesh gap by an independent ordinary Brownian bridge, and let
\(X^{(k)}=(X_s^{(k)})_{0\le s\le b_{T,k}}\) denote the resulting continuous
block path. Thus
\[
  X_{rh_T}^{(k)}=Y_r,
  \qquad
  r=0,\ldots,n_{T,k}.
\]
Under the base candidate law, the completed process \(X^{(k)}\) is an
ordinary Brownian motion started from \(Y_0\) and run for time \(b_{T,k}\).
Set $a_{\kappa,\zeta}
  :=
  [1+2(\kappa+\zeta)]^{-d},$
and define the good mesh event
\[
  \Omega_{T,k}^{\mathrm{mesh}}(\kappa)
  :=
  \left\{
    \max_{1\le r\le n_{T,k}}
    \sup_{(r-1)h_T\le s\le rh_T}
    \left\|
      X_s^{(k)}-Y_r
    \right\|_\infty
    \le \kappa\delta_T
  \right\}.
\]
Then the following statements hold uniformly over the incoming history and
the phase of the \(\delta_T\)-grid.
\begin{itemize}
    \item (1) on \(\Omega_{T,k}^{\mathrm{mesh}}(\kappa)\),
\begin{equation}
  U_{T,k}^{\raw}
  \ge
  a_{\kappa,\zeta}
  \left|
    W_{[0,b_{T,k}]}^{\zeta\delta_T}
      (X^{(k)})
  \right|.
  \label{eq:raster-containment}
\end{equation}
\item (2) for a dimension-dependent constant \(C_d<\infty\),
\begin{equation}
  \Pp\!\left(
    \Omega_{T,k}^{\mathrm{mesh}}(\kappa)^c
    \,\middle|\,
    \Xi_k
  \right)
  \le
  4dN_T
  \exp\left\{
    -\frac{\kappa^2\delta_T^2}{C_dh_T}
  \right\}.
  \label{eq:app-mesh-complement}
\end{equation}
\item (3) for \(v>0\), define $s_{T,k}
  :=
  b_{T,k}\delta_T^{-2}$ and $f_{T,k}(v)
  :=
  \frac{v}{a_{\kappa,\zeta}\delta_T^d}.$ 
Then Brownian scaling and \eqref{eq:raster-containment} give
\begin{equation}
  \begin{split}
  \Pp\{U_{T,k}^{\raw}\le v\mid\Xi_k\}
  \le
  \Pp\!\left(
    \Omega_{T,k}^{\mathrm{mesh}}(\kappa)^c
    \,\middle|\,
    \Xi_k
  \right)+
  \Pp\left\{
    |W_{s_{T,k}}^\zeta|
    \le f_{T,k}(v)
  \right\},
  \end{split}
  \label{eq:app-raw-to-source}
\end{equation}
and the corresponding fixed-radius exponent satisfies
\begin{equation}
  \frac{\lambda_d}{2}
  s_{T,k}f_{T,k}(v)^{-2/d}
  =
  b_{T,k}
  a_{\kappa,\zeta}^{2/d}
  \Lambda_{\FK}(v).
  \label{eq:app-source-exponent}
\end{equation}
\end{itemize}
\end{lemma}

\begin{proof}
The Gaussian endpoint kernel gives independent increments
\[
  Y_r-Y_{r-1}
  \sim N(0,h_TI_d),
  \qquad
  r=1,\ldots,n_{T,k}.
\]
Completing each gap by an independent Brownian bridge conditional on its two
endpoints is the standard disintegration of Brownian motion into its mesh
values and conditional bridges. Consequently, after the bridges are added,
\(X^{(k)}\) is an ordinary Brownian motion started from \(Y_0\). Integrating
out the bridges recovers the original product-Gaussian endpoint kernel.

We first prove \eqref{eq:raster-containment}. Let $z\in
  W_{[0,b_{T,k}]}^{\zeta\delta_T}(X^{(k)}).$
Then there exists \(s\in[0,b_{T,k}]\) such that $|z-X_s^{(k)}|_2\le\zeta\delta_T.$
Choose \(r\in\{1,\ldots,n_{T,k}\}\) so that $(r-1)h_T\le s\le rh_T.$
On \(\Omega_{T,k}^{\mathrm{mesh}}(\kappa)\), $\|X_s^{(k)}-Y_r\|_\infty
  \le\kappa\delta_T.$
Since Euclidean distance bounds coordinatewise distance, $\|z-X_s^{(k)}\|_\infty
  \le\zeta\delta_T,$
and hence
\[
  \|z-Y_r\|_\infty
  \le(\kappa+\zeta)\delta_T.
\]

Let $m_r=\left\lfloor Y_r/\delta_T\right\rfloor.$
Since $Y_r\in
  \delta_T(m_r+[0,1)^d),$
the point \(z\) belongs to the coordinatewise enlargement $\delta_T
  \bigl(
    m_r+[-(\kappa+\zeta),\,1+\kappa+\zeta]^d
  \bigr).$
This enlarged cell has volume $[1+2(\kappa+\zeta)]^d\delta_T^d.$
Taking the union over the distinct grid cells visited by
\(Y_1,\ldots,Y_{n_{T,k}}\), we obtain
\[
  \left|
    W_{[0,b_{T,k}]}^{\zeta\delta_T}(X^{(k)})
  \right|
  \le
  [1+2(\kappa+\zeta)]^d
  \delta_T^d
  \left|
    \left\{
      \left\lfloor Y_r/\delta_T\right\rfloor:
      1\le r\le n_{T,k}
    \right\}
  \right|.
\]
The last factor is
\(U_{T,k}^{\raw}/\delta_T^d\), which proves
\eqref{eq:raster-containment}.

We next prove \eqref{eq:app-mesh-complement}. Fix one coordinate and one
mesh gap. Under the completed Brownian law, time reversal on that gap shows
that $\left(
    X_{rh_T}^{(k)}
    -
    X_{rh_T-s}^{(k)}
  \right)_{0\le s\le h_T}$
has the law of a one-dimensional Brownian motion on \([0,h_T]\).
The reflection principle and a Gaussian tail bound therefore give
\[
  \Pp\left\{
    \sup_{0\le s\le h_T}|B_s|
    >\kappa\delta_T
  \right\}
  \le
  4\exp\left\{
    -\frac{\kappa^2\delta_T^2}{C_dh_T}
  \right\}.
\]
There are at most \(d\,n_{T,k}\le dN_T\) coordinate-gap pairs. A union
bound over them proves \eqref{eq:app-mesh-complement}. The estimate is
uniform in the incoming endpoint because Brownian increments are
translation invariant.

It remains to prove the source translation. On
\(\Omega_{T,k}^{\mathrm{mesh}}(\kappa)\), equation \eqref{eq:raster-containment} implies
\[
  \{U_{T,k}^{\raw}\le v\}
  \subset
  \left\{
    \left|
      W_{[0,b_{T,k}]}^{\zeta\delta_T}(X^{(k)})
    \right|
    \le
    \frac{v}{a_{\kappa,\zeta}}
  \right\}.
\]
Combining this inclusion with the complement of the good event gives
\[
  \begin{aligned}
  \Pp\{U_{T,k}^{\raw}\le v\mid\Xi_k\}
  \le
  \Pp\!\left(
    \Omega_{T,k}^{\mathrm{mesh}}(\kappa)^c
    \,\middle|\,
    \Xi_k
  \right)+
  \Pp\left\{
    \left|
      W_{[0,b_{T,k}]}^{\zeta\delta_T}
    \right|
    \le
    \frac{v}{a_{\kappa,\zeta}}
  \right\}.
  \end{aligned}
\]

Translate the completed Brownian path to start from the origin and
rescale space by the factor \(\delta_T^{-1}\). The rescaled process
(here \(\widehat B\) is different from the process in Section~\ref{sec:setting})
\[
  \widehat B_u
  :=
  \delta_T^{-1}
  \left(
    X_{\delta_T^2u}^{(k)}-Y_0
  \right),
  \qquad
  0\le u\le s_{T,k},
\]
is standard Brownian motion, and
\[
  \left|
    W_{[0,b_{T,k}]}^{\zeta\delta_T}(X^{(k)})
  \right|
  =
  \delta_T^d
  |W_{s_{T,k}}^\zeta(\widehat B)|.
\]
Thus
\[
  \left\{
    \left|
      W_{[0,b_{T,k}]}^{\zeta\delta_T}
    \right|
    \le
    \frac{v}{a_{\kappa,\zeta}}
  \right\}
  =
  \left\{
    |W_{s_{T,k}}^\zeta|
    \le
    \frac{v}{a_{\kappa,\zeta}\delta_T^d}
  \right\},
\]
which proves \eqref{eq:app-raw-to-source}.

Finally,
\[
  \begin{aligned}
  \frac{\lambda_d}{2}
  s_{T,k}f_{T,k}(v)^{-2/d}
  =
  \frac{\lambda_d}{2}
  \left(
    b_{T,k}\delta_T^{-2}
  \right)
  \left(
    \frac{v}
         {a_{\kappa,\zeta}\delta_T^d}
  \right)^{-2/d}=
  b_{T,k}
  a_{\kappa,\zeta}^{2/d}
  \frac{\lambda_d}{2}v^{-2/d}=
  b_{T,k}
  a_{\kappa,\zeta}^{2/d}
  \Lambda_{\FK}(v),
  \end{aligned}
\]
which is \eqref{eq:app-source-exponent}.
\end{proof}

\begin{lemma}[Freeze--then--wait uniformization and error ledger]
\label{lem:freeze-wait-uniformization}
There exist deterministic step functions
\[
  u_T,\theta_T,\kappa_T,\zeta_T\downarrow0,
  \qquad
  A_T^{\cut},v_T^{\up}\uparrow\infty,
  \qquad
  v_T^{\up}=A_T^{\cut}+J+3,
\]
and deterministic nonnegative remainders \(r_{T,k}\) such that, uniformly
over bulk blocks \(q_{k+1}\ge\theta_T\), incoming histories, grid phases,
and \(v\in[u_T,v_T^{\up}]\),
\[
  \log\Pp\{U_{T,k}^{\raw}\le v\mid\Xi_k\}
  \le
  -b_{T,k}\widetilde\lambda_T(v)+r_{T,k}.
\]
Moreover,
\[
  u_T\le\frac14v_\star\theta_T^p,
  \qquad
  \sup_{k<K_T}\frac{r_{T,k}}{b_T}\longrightarrow0,
  \qquad
  \sum_{k<K_T}r_{T,k}=o(T),
\]
and the same choice of auxiliary sequences satisfies the error ledger
\eqref{eq:main-ledger}.
In addition, the deterministic candidate multiplicities satisfy $\log\left(\sum_{k<K_T}M_{T,k}\right)=o(T).$
\end{lemma}

\begin{proof}
The proof is mainly a diagonal argument. We construct the auxiliary quantities by freezing all parameters at an
integer stage and then waiting until the corresponding uniform estimates
and error bounds hold.

\medskip
\noindent
\textit{Step 1: frozen parameters.}
Choose decreasing positive sequences $\theta_n\downarrow0$ and $u_n\downarrow0$
such that $\theta_n^p\le n^{-3}$ and $ u_n
  \le
  \min\left\{
    n^{-3},
    \frac14v_\star\theta_n^p
  \right\}.$
Choose \(A_n\uparrow\infty\) so that $\Lambda_{\FK}(A_n)\le n^{-3},$
and put $R_n:=A_n+J+3.$
Finally, choose decreasing positive sequences
\(\kappa_n,\zeta_n\downarrow0\), with
\(\kappa_n,\zeta_n\le n^{-1}\), such that $(\kappa_n+\zeta_n)\Lambda_{\FK}(u_n)
  \le n^{-3}.$
Set $a_n
  :=
  [1+2(\kappa_n+\zeta_n)]^{-d}.$
Thus, at stage \(n\), the threshold interval
\([u_n,R_n]\), the comparison radius \(\zeta_n\), the mesh tolerance
\(\kappa_n\), and the dilation factor \(a_n\) are all fixed before the
fixed-radius source theorem is invoked.

\medskip
\noindent
\textit{Step 2: uniformity at one frozen stage.}
Fix \(n\). The normalized balanced durations \(b_{T,k}/b_T\) remain,
for all sufficiently large \(T\), in a compact subset of
\((0,\infty)\), and the threshold \(v\) lies in the compact interval
\([u_n,R_n]\). The comparison radius \(\zeta_n>0\) remains fixed.

We claim that, beyond a deterministic stage-dependent horizon,
the fixed-radius source estimate can be made uniform over these
parameters. More precisely, with $s_{T,k}
  =
  b_{T,k}\delta_T^{-2}$ and $f_{T,k}^{(n)}(v)
  =
  \frac{v}{a_n\delta_T^d},$
one has
\begin{equation}
  \log\Pp\left\{
    |W_{s_{T,k}}^{\zeta_n}|
    \le f_{T,k}^{(n)}(v)
  \right\}
  \le
  -b_{T,k}a_n^{2/d}\Lambda_{\FK}(v)
  +n^{-2}b_T
  \label{eq:app-frozen-source}
\end{equation}
uniformly over balanced block durations and
\(v\in[u_n,R_n]\).

Indeed, if no such uniform horizon existed, there would be a sequence
\(T_j\to\infty\), balanced blocks \(k_j\), and thresholds
\(v_j\in[u_n,R_n]\) for which
\eqref{eq:app-frozen-source} failed. Let $s_j:=b_{T_j,k_j}\delta_{T_j}^{-2}$ and $F_j:=\frac{v_j}{a_n\delta_{T_j}^d}.$
Since \(s_j\to\infty\), pass to a subsequence on which \(s_j\) is
strictly increasing and \(s_j>e^2\). For \(s\ge e^2\), set
\[
  g_d(s):=
  \begin{cases}
    s/\log s,&d=2,\\
    s,&d\ge3.
  \end{cases}
\]
The schedule estimates at the start of this appendix, the bounds
\(u_n\le v_j\le R_n\), and the uniform comparability
\(b_{T_j,k_j}\asymp b_{T_j}\) imply $F_j\longrightarrow\infty$ and $\frac{F_j}{g_d(s_j)}\longrightarrow0.$
Define one deterministic function \(f\) by \(f(s_j)=F_j\), by
\(f(s)=\sqrt{g_d(s)}\) at all other \(s\ge e^2\), and by
\(f(s)=1\) for \(s<e^2\). Then \(f(s)\to\infty\) and
\(f(s)=o(g_d(s))\). The fixed-radius theorem
\eqref{eq:hard-lower-tail} applies to this function; its upper bound
also applies to the closed-ball sausage, as explained in
Section~\ref{sec:efficiency}. Along the selected sequence it yields
\[
  \begin{aligned}
  \log\Pp\{|W_{s_j}^{\zeta_n}|\le F_j\}
  &\le
  -\frac{\lambda_d}{2}s_jF_j^{-2/d}
  +o(s_jF_j^{-2/d})=
  -b_{T_j,k_j}a_n^{2/d}\Lambda_{\FK}(v_j)
  +o(b_{T_j}).
  \end{aligned}
\]
Here the last error bound uses \(v_j\ge u_n>0\) and the bounded
ratios \(b_{T_j,k_j}/b_{T_j}\). This contradicts the failure of
\eqref{eq:app-frozen-source}, whose allowed error is the fixed positive
multiple \(n^{-2}b_{T_j}\). 

Lemma~\ref{lem:raster-source-translation} is uniform in the incoming
endpoint and grid phase, so the same frozen estimate applies after the
grid-to-source translation.

\medskip
\noindent
\textit{Step 3: waiting horizons.}
Choose deterministic horizons $T_1<T_2<\cdots$ and $T_n\uparrow\infty$
large enough that, for every \(T\ge T_n\), the frozen estimate
\eqref{eq:app-frozen-source} holds and all of the following conditions are
satisfied.
First, we require the scaled lower threshold to diverge: $\frac{u_n}{a_n\delta_T^d}\ge n.$
Next require the spectral-certification bounds
\begin{equation*}
  |\chic_{d,T}-\chi_d|u_n^{-2/d}
  \le n^{-3},
  \qquad
  \widetilde\lambda_T(A_n)
  \le2n^{-3}.
\end{equation*}

The uniform source-range conditions over
\(v\in[u_n,R_n]\) are enforced by
\begin{equation*}
  \frac{R_n}
       {a_n b_T\delta_T^{d-2}}
  \le n^{-1},
  \qquad d\ge3,
\end{equation*}
and
\begin{equation*}
  \frac{
    R_n\log(b_{T,\max}\delta_T^{-2})
  }{
    a_n b_T
  }
  \le n^{-1},
  \qquad d=2.
\end{equation*}
These are precisely the conditions obtained by substituting
\[
  s_{T,k}=b_{T,k}\delta_T^{-2},
  \qquad
  f_{T,k}^{(n)}(v)
  =\frac{v}{a_n\delta_T^d}
\]
into the two dimension-dependent source hypotheses, using the uniform
comparability \(b_{T,k}\asymp b_T\).
Next, enlarge \(T_n\) further so that the mesh-complement bound from
Lemma~\ref{lem:raster-source-translation} satisfies
\begin{equation}
  \frac{\kappa_n^2\delta_T^2}{C_dh_T}
  \ge
  2b_{T,\max}\Lambda_{\FK}(u_n)
  +nb_T
  +\log(4dN_T).
  \label{eq:app-mesh-domination}
\end{equation}
Consequently, uniformly over balanced blocks and incoming histories,
\[
  \begin{aligned}
  \Pp(\Omega_{T,k}^{\mathrm{mesh}}(\kappa_n)^c\mid\Xi_k)
  &\le E_{T,n},
  \end{aligned}
\]where $E_{T,n}:=
  \exp\{-2b_{T,\max}\Lambda_{\FK}(u_n)-nb_T\}.$
For \(v\in[u_n,R_n]\), denote the exponential upper bound from
\eqref{eq:app-frozen-source} by
\[
  B_{T,k,n}(v):=
  \exp\{-b_{T,k}a_n^{2/d}\Lambda_{\FK}(v)+n^{-2}b_T\}.
\]
Since \(a_n\le1\), \(b_{T,k}\le b_{T,\max}\), and
\(\Lambda_{\FK}(v)\le\Lambda_{\FK}(u_n)\),
\[
  \frac{E_{T,n}}{B_{T,k,n}(v)}
  \le
  \exp\{-b_{T,\max}\Lambda_{\FK}(u_n)-(n+n^{-2})b_T\}
  \le e^{-nb_T}.
\]
Thus the mesh error is exponentially smaller than the source exponential
upper bound, uniformly on the frozen interval.

Finally, enlarge \(T_n\) again so that, for every \(T\ge T_n\), each of
\[
  \frac{K_T\log(N_T+3)}{T},
  \qquad
  \widetilde\lambda_T(A_n),
  \qquad
  \theta_n^p,
  \qquad
  K_T^{-p},
  \qquad
  \frac{K_T}{T},
  \qquad
  \gamma_T^2,
  \qquad
  |\Jc_T-J|
\]
is at most \(n^{-2}\).

\medskip
\noindent
\textit{Step 4: the freeze--then--wait diagonal.}
For \(T\ge T_1\), define
\[
  n(T):=\max\{n:T_n\le T\};
\]
for \(T<T_1\), use stage \(1\). Set
\[
  \begin{aligned}
  u_T&:=u_{n(T)},&
  \theta_T&:=\theta_{n(T)},&
  \kappa_T&:=\kappa_{n(T)},&
  \zeta_T&:=\zeta_{n(T)},\\
  A_T^{\cut}&:=A_{n(T)},&
  v_T^{\up}&:=R_{n(T)}.
  \end{aligned}
\]
Since \(n(T)\to\infty\), these deterministic step functions have all the
limits claimed in the lemma, and
\[
  v_T^{\up}
  =
  A_T^{\cut}+J+3.
\]

Fix \(T\ge T_1\), write \(n=n(T)\), and take
\(v\in[u_n,R_n]\). Lemma~\ref{lem:raster-source-translation} gives
\[
  \begin{aligned}
  \Pp\{U_{T,k}^{\raw}\le v\mid\Xi_k\}
  \le
  \Pp(\Omega_{T,k}^{\mathrm{mesh}}(\kappa_n)^c\mid\Xi_k)+
  \Pp\left\{
    |W_{s_{T,k}}^{\zeta_n}|
    \le f_{T,k}^{(n)}(v)
  \right\}.
  \end{aligned}
\]
The two terms are bounded by \(E_{T,n}\) and \(B_{T,k,n}(v)\),
respectively. The comparison established after
\eqref{eq:app-mesh-domination} therefore gives
\[
  \Pp\{U_{T,k}^{\raw}\le v\mid\Xi_k\}
  \le B_{T,k,n}(v)(1+e^{-nb_T})
  \le 2B_{T,k,n}(v).
\]
Moreover,
\[
  1-a_n^{2/d}
  =
  1-[1+2(\kappa_n+\zeta_n)]^{-2}
  \le C_d(\kappa_n+\zeta_n),
\]
so, since \(v\ge u_n\) and
\(\Lambda_{\FK}\) is decreasing,
\[
  \bigl(1-a_n^{2/d}\bigr)\Lambda_{\FK}(v)
  \le
  C_d(\kappa_n+\zeta_n)\Lambda_{\FK}(u_n)
  \le C_dn^{-3}.
\]
The spectral-certification condition gives
\[
  |\widetilde\lambda_T(v)-\Lambda_{\FK}(v)|
  \le
  |\chic_{d,T}-\chi_d|u_n^{-2/d}
  \le n^{-3}.
\]
Hence the source exponential upper bound satisfies
\[
  B_{T,k,n}(v)\le\exp\left\{
    -b_{T,k}\widetilde\lambda_T(v)
    +C_d(n^{-2}+n^{-3})b_T
  \right\}.
\]
Combining the last two bounds introduces at most the additive term
\(\log2\) after taking logarithms. We may therefore choose nonnegative
\(r_{T,k}\) such that
\[
  \log\Pp\{U_{T,k}^{\raw}\le v\mid\Xi_k\}
  \le
  -b_{T,k}\widetilde\lambda_T(v)+r_{T,k}
\]
and
\[
  \sup_{k<K_T}\frac{r_{T,k}}{b_T}
  \le
  C_d(n(T)^{-2}+n(T)^{-3})
  +\frac{\log2}{b_T}
  \longrightarrow0.
\]
On the bounded range \(T<T_1\), enlarge the remainders deterministically
if necessary; this does not affect any asymptotic conclusion.

Since \(b_{T,k}\asymp b_T\) uniformly and $\sum_{k<K_T}b_{T,k}=T,$
we have \(K_Tb_T=O(T)\). Therefore
\[
  \sum_{k<K_T}r_{T,k}
  \le
  K_Tb_T
  \sup_{k<K_T}\frac{r_{T,k}}{b_T}
  =
  o(T).
\]

The final waiting conditions imply
\[
  \begin{aligned}
  K_T\log(N_T+3)
  +T\widetilde\lambda_T(A_T^{\cut})
  +T\theta_T^p
  +TK_T^{-p}
  +K_T
  +T\gamma_T^2
  +T|\Jc_T-J|
  =o(T),
  \end{aligned}
\]
which is exactly \eqref{eq:main-ledger}.

\medskip
\noindent
\textit{Step 5: deterministic candidate multiplicity.}
It remains to record the subexponential multiplicity estimate. Concavity of
\(q\mapsto q^p\) and the balanced mesh give
\[
  \max_{k<K_T}\calL_{T,k}
  \le
  C_dT(K_T^{-1}+N_T^{-1})^p
  =
  o(T).
\]
Also,
\[
  \log C_{T,d}^{\diag}
  =
  O_d(1)
  +
  O_d(\gamma_T^2b_{T,\max})
  =
  o(T).
\]
Since $X_{T,k}
  =
  C_{T,d}^{\diag}K_T^2e^{2\calL_{T,k}}
  \ge2$
and \(M_{T,k}=\lceil X_{T,k}\rceil\),
\[
  M_{T,k}
  <
  X_{T,k}+1
  \le
  \frac32X_{T,k}.
\]
Consequently,
\[
  \sum_{k<K_T}M_{T,k}
  \le
  \frac32
  C_{T,d}^{\diag}
  K_T^3
  e^{2\max_{k<K_T}\calL_{T,k}}.
\]
Taking logarithms and using
\(\log K_T=o(T)\) proves
\[
  \log\left(\sum_{k<K_T}M_{T,k}\right)=o(T).
\]
\end{proof}

\begin{lemma}[Uniform continuation estimates]
\label{lem:continuation-estimates}
Use the bulk and terminal block sets and the quantities
\(\eta_{T,k}^{\mathrm{bulk}}\) and \(\eta_{T,k}^{\mathrm{term}}\) defined in
\eqref{eq:bulk-eta}--\eqref{eq:terminal-eta}. For the auxiliary sequences and
remainders supplied by Proposition~\ref{prop:fresh-tail}, the following
statements hold for all sufficiently large \(T\).
For every \(k\in\mathcal B_T\),
\begin{equation}
  b_{T,k}\widetilde\lambda_T(u_T)
  \ge
  3^{2/d}\calL_{T,k}
  >
  \calL_{T,k}.
  \label{eq:app-release-margin}
\end{equation}

Let \(\mathcal D_{T,k}\) be the event that at least one endpoint in the
fresh candidate block receives the dagger token. On
\(\mathcal D_{T,k}^c\),
\begin{equation}
  V\ge U_{T,k}^{\raw},
  \qquad
  \{V<v\}
  \subset
  \{U_{T,k}^{\raw}<v\}.
  \label{eq:history-removal}
\end{equation}

Moreover, uniformly over every reachable history,
\begin{equation}
  \E\!\left[
    W_{T,k}^{\esc}
    \,\middle|\,
    \Xi_k
  \right]
  +
  \E\!\left[
    W_{T,k}^{\mathrm{vol}};
    \mathcal D_{T,k}
    \,\middle|\,
    \Xi_k
  \right]
  \le
  C_d
  \exp\{C_d\gamma_T^2b_{T,k}\}.
  \label{eq:cemetery-bounds}
\end{equation}

Moreover,
\begin{equation}
  \sum_{k\in\mathcal B_T}
  \eta_{T,k}^{\mathrm{bulk}}
  +
  \sum_{k\in\mathcal T_T}
  \eta_{T,k}^{\mathrm{term}}
  =
  o(T).
  \label{eq:app-continuation-sum}
\end{equation}
\end{lemma}

\begin{proof}
We prove the release margin, the history-removal relation, the
escape/window estimate, and finally the accumulated remainder bound.

\medskip
\noindent
\textit{Step 1: release margin on bulk blocks.}
Fix \(k\in\mathcal B_T\), and write $q=q_k$ and $q'=q_{k+1}.$
Since \(q'\ge\theta_T\), concavity of \(q\mapsto q^p\) gives
\[
  q^p-(q')^p
  \le
  p(q')^{p-1}(q-q').
\]
Using \(T(q-q')=b_{T,k}\), we obtain
\begin{equation}
  \calL_{T,k}
  =
  T\Jc_T\bigl(q^p-(q')^p\bigr)
  \le
  b_{T,k}p\Jc_T(q')^{p-1}.
  \label{eq:app-release-upper}
\end{equation}

The exported coupling in \eqref{eq:fresh-tail-coupling} gives
\[
  u_T
  \le
  \frac14v_\star\theta_T^p
  \le
  \frac14v_\star(q')^p.
\]
Because \(\widetilde\lambda_T(v)=\chic_{d,T}v^{-2/d}\),
\[
  \widetilde\lambda_T(u_T)
  \ge
  \chic_{d,T}
  \left(
    \frac14v_\star(q')^p
  \right)^{-2/d}.
\]
For the exact constants,
\[
  \chi_dv_\star^{-2/d}
  =
  \frac d2v_\star
  =
  pJ,
  \qquad
  -\frac{2p}{d}=p-1.
\]
Consequently, the exact release-to-tail ratio is \(4^{2/d}\):
\[
  \chi_d
  \left(
    \frac14v_\star(q')^p
  \right)^{-2/d}
  =
  4^{2/d}pJ(q')^{p-1}.
\]
Since $\chic_{d,T}\longrightarrow\chi_d$ and $\Jc_T\longrightarrow J,$
the fixed slack \(3^{2/d}<4^{2/d}\) implies, after increasing the
deterministic waiting horizon if necessary,
\[
  \widetilde\lambda_T(u_T)
  \ge
  3^{2/d}p\Jc_T(q')^{p-1}.
\]
Combining this inequality with \eqref{eq:app-release-upper} proves
\eqref{eq:app-release-margin}.

\medskip
\noindent
\textit{Step 2: removal of previously paid tokens.}
For a fresh candidate \(Y\), let
\[
 \mathcal E_{T,k}^{\raw}(Y)
  :=
  \left\{
    \left\lfloor Y_r/\delta_T\right\rfloor:
    1\le r\le n_{T,k}
  \right\}
\]
be its set of distinct infinite-grid cells. On
\(\mathcal D_{T,k}^c\), every raw cell is valid, so $\mathcal E_T(Y)
  =
  \mathcal E_{T,k}^{\raw}(Y).$
The updated paid-token set is therefore
\[
  D_k^Y
  =
  D_k
  \cup
  \bigl(
    \mathcal E_{T,k}^{\raw}(Y)\setminus D_k
  \bigr)
  =
  D_k\cup\mathcal E_{T,k}^{\raw}(Y).
\]
Hence 
\[
  \begin{aligned}
  V
  =
  A+w=
  \delta_T^d|D_k|
  +
  \delta_T^d
  \left|
    \mathcal E_{T,k}^{\raw}(Y)\setminus D_k
  \right|=
  \delta_T^d
  \left|
    D_k\cup\mathcal E_{T,k}^{\raw}(Y)
  \right|\ge
  \delta_T^d
  \left|
    \mathcal E_{T,k}^{\raw}(Y)
  \right|
  =
  U_{T,k}^{\raw}.
  \end{aligned}
\]
The event inclusion in \eqref{eq:history-removal} follows immediately.

\medskip
\noindent
\textit{Step 3: escape and window-exit terms.}
Fix the current endpoint \(x=x(\Xi_k)\), and let $Z\sim N(0,b_{T,k}I_d)$
be the total endpoint displacement of the fresh candidate. Since
\(d_T^{\cem}\) is \(1\)-Lipschitz in the \(\ell^\infty\) norm,
\[
  d_T^{\cem}(x)-d_T^{\cem}(x+Z)
  \le
  \|Z\|_\infty.
\]
Therefore
\[
  \frac{
    H_{T,k+1}^{\esc}(x+Z)
  }{
    H_{T,k}^{\esc}(x)
  }
  \le
  e^{\gamma_T\|Z\|_\infty}.
\]
Because \(g_{T,k}\le1\) and
\(H_{T,k}(\Xi_k)\ge H_{T,k}^{\esc}(x)\), $W_{T,k}^{\esc}
  \le
  e^{\gamma_T\|Z\|_\infty}.$
A union bound over coordinates and Gaussian square completion gives
\[
  \E e^{\gamma_T\|Z\|_\infty}
  \le
  C_d
  e^{C_d\gamma_T^2b_{T,k}}.
\]

It remains to control the volume branch on
\(\mathcal D_{T,k}\). Put $r=d_T^{\cem}(x).$
Since \(g_{T,k}\le1\),
\(H_{T,k+1}^{\mathrm{vol}}\le1\), and
\(H_{T,k}(\Xi_k)\ge H_{T,k}^{\esc}(x)=e^{-\gamma_Tr}\), then $W_{T,k}^{\mathrm{vol}}
  \le
  e^{\gamma_Tr}.$
If \(x\) is inside the inward window, the event
\(\mathcal D_{T,k}\) requires an \(\ell^\infty\) displacement of at least
\(r\) during the block. A reflection estimate and a Gaussian tail bound
therefore give
\[
  \Pp(\mathcal D_{T,k}\mid\Xi_k)
  \le
  C_d
  \exp\left\{
    -\frac{r^2}{C_db_{T,k}}
  \right\}.
\]
Consequently,
\[
  \begin{aligned}
  \E\!\left[
    W_{T,k}^{\mathrm{vol}};
    \mathcal D_{T,k}
    \,\middle|\,
    \Xi_k
  \right]
  &\le
  C_d
  \exp\left\{
    \gamma_Tr
    -\frac{r^2}{C_db_{T,k}}
  \right\}\le
  C_d
  e^{C_d\gamma_T^2b_{T,k}},
  \end{aligned}
\]
where the last step follows by completing the square in \(r\). If \(x\)
is outside the inward window, then \(r=0\), and the same bound is immediate.
Combining the two estimates proves \eqref{eq:cemetery-bounds}.

\medskip
\noindent
\textit{Step 4: accumulated bulk and terminal remainders.}
For the bulk blocks, combining \eqref{eq:main-ledger} with $\sum_{k<K_T}b_{T,k}=T$ 
gives
\[
  \begin{aligned}
  \sum_{k\in\mathcal B_T}
  \eta_{T,k}^{\mathrm{bulk}}
  &\le
  \sum_{k<K_T}r_{T,k}
  +
  K_T\log(N_T+3)
  +
  T\widetilde\lambda_T(A_T^{\cut})+
  C_dT\gamma_T^2
  +
  C_dK_T
  =
  o(T).
  \end{aligned}
\]

For the terminal blocks, suppose first that
\(\mathcal T_T\ne\varnothing\), and let \(k_0\) be its first index.
The balanced schedule gives $q_{k_0}
  \le
  \theta_T+\frac{C_d}{K_T}.$
Since the releases telescope, $\sum_{k\in\mathcal T_T}\calL_{T,k}
  =
  T\Jc_Tq_{k_0}^p\le
  T\Jc_T
  \left(
    \theta_T^p+C_dK_T^{-p}
  \right).$

Also, from the definition of \(C_{T,d}^{\diag}\),
\[
  \log C_{T,d}^{\diag}
  \le
  C_d+C_d\gamma_T^2b_{T,\max}.
\]
The balanced block lengths satisfy $K_Tb_{T,\max}=O(T),$
and therefore
\[
  \begin{aligned}
  \sum_{k\in\mathcal T_T}
  \frac12\log C_{T,d}^{\diag}
\le
  C_dK_T
  +
  C_d\gamma_T^2K_Tb_{T,\max}
  &\le
  C_dK_T+C_dT\gamma_T^2
  =
  o(T).
  \end{aligned}
\]
Together with \eqref{eq:main-ledger}, this proves $\sum_{k\in\mathcal T_T}
  \eta_{T,k}^{\mathrm{term}}
  =
  o(T).$
If \(\mathcal T_T=\varnothing\), the terminal sum is zero. Combining the
bulk and terminal estimates proves \eqref{eq:app-continuation-sum}.

All statements above are required only after a deterministic threshold
\(T_0(d,a)\). On the bounded range of smaller values of \(T\), the
remainders may be enlarged deterministically without affecting any
asymptotic conclusion. The finite-\(T\) estimator law itself is unchanged.
\end{proof}

\section{Additional details of the lattice surrogate experiments}
\label{sec:numericsdetail}

The complete implementation, fixed configurations, and instructions for reproducing all experiments are publicly available at the Github repository. 


\subsection{Natural-radius scaling and lattice surrogate}

The experiments take the natural-radius multiplier to be $a=1$. Throughout
Section~\ref{sec:numerics}, $T$ denotes the natural-radius time horizon
in the natural-radius formulation, and $\eps_T=T^{-1/d}$. Choose the
lattice spacing $h=\eps_T/5$. A displacement by one lattice site represents
a physical displacement of size $h$, and a lattice ball of radius five
represents the physical radius $5h=\eps_T$.

Let $S=(S_s)_{s\geq0}$ be a rate-one continuous-time Markov
chain on $\mathbb Z^d$, started from the origin, with uniform jump
directions. Since $\operatorname{Cov}(S_s)=(s/d)I_d$, the rescaled process
$X_t:=hS_{dt/h^2}$ satisfies $\operatorname{Cov}(X_t)=tI_d$, the covariance
normalization of the Brownian motion in the continuous model. The natural-radius time horizon $T$ therefore corresponds to the raw lattice horizon $t_{\mathrm{raw}}=\frac{dT}{h^2}=25d\,T^{1+2/d}.$
The cases reported in Section~\ref{sec:numerics} are indexed by the natural-radius time horizon $T$. In raw lattice coordinates, the shrinking natural radius has
become the fixed discrete radius five.

Let $\mathcal B_5:=\{z\in\mathbb Z^d:\lVert z\rVert_2\leq5\}$, and
define the discrete sausage
$A:=\bigcup_{0\leq s\leq t_{\mathrm{raw}}}(S_s+\mathcal B_5)$.
Each occupied lattice site represents physical volume $h^d$, so the
continuous penalty $T|W_T^{\eps_T}|$ is represented by $Th^d|A|$.
Setting $\beta:=Th^d=5^{-d}$ gives the lattice target $\alpha^{\mathrm{lat}}_{d,T}
 :=\mathbb E_{\mathrm{lat}}\!\left[\exp\{-\beta|A|\}\right],$
where $\mathbb E_{\mathrm{lat}}$ denotes expectation under the rate-one
nearest-neighbor chain. The visited-site sequence and the occupied-site
union $A$ are computed exactly for this finite lattice model.

\subsection{Relation to the continuous construction}

The lattice estimator uses candidate augmentation as in
Algorithm~\ref{alg:package-w}. At block $k$, it generates $M_k$ independent
base-kernel candidate continuations, evaluates their exact new occupied
volumes and positive continuation guides, and selects one candidate. Let
$\pi_{k,i}>0$ be candidate $i$'s actual implemented selection probability,
including any defensive mixing or integer-ticket rounding. If $I_k$ is the
selected label, the retained likelihood correction is $L_k=\frac{1}{M_k\pi_{k,I_k}}.$
Conditional on the candidate cloud, this correction turns the selected
candidate's contribution into the uniform average over the cloud.
Iterating over blocks shows that the final payoff
$\exp\{-\beta|A|\}\prod_k L_k$ is unbiased for
$\alpha^{\mathrm{lat}}_{d,T}$, independently of the guide approximation.

Several components of the continuous construction are unnecessary for the
lattice surrogate. The occupied set is the exact discrete sausage of the
selected lattice path, so no inward grid approximation or valid-token
selection is required. The jump count and jump directions determine the
visited-site sequence and hence the entire discrete sausage, so no
Brownian-bridge completion is needed. Finally, the terminal payoff uses the
exact occupied set $A$, leaving no uncharged residual region and requiring
no residual Poisson-void correction. The experiments thus study an exact
finite-horizon estimator of the lattice target.

\section*{Acknowledgments}
The authors used artificial intelligence tools for language editing and
coding. They also used
PaperFormalizer~\cite{TanBlanchetSyrgkanis2026PaperFormalizer}
for assisted formalization and proof verification in Lean.
The formalized statements corresponding to all results examined
passed Lean verification.
The authors take full responsibility for all mathematical arguments,
results, and conclusions presented in the paper.
J. Blanchet gratefully acknowledges support from ONR under award
N00014-24-1-2655, and the National Science Foundation (NSF) under grants
2312204 and 2403007.

\end{document}